\documentclass[a4paper]{amsart}
\usepackage[utf8]{inputenc}
\usepackage[english]{babel}

\usepackage{amscd,amssymb}
\usepackage{amsmath}
\usepackage{amsthm}
\usepackage{graphicx}
\usepackage{fancybox}
\usepackage{epic,eepic}
\usepackage{amstext}
\usepackage{mathtools}
\usepackage{esint}
\usepackage{booktabs}
\usepackage[hide links]{hyperref}
\usepackage{comment}
\usepackage{mathtools}
\usepackage{tikz,pgfplots}
\usetikzlibrary{automata, positioning, arrows, decorations.pathreplacing, calc, chains, shapes.misc}
\usepackage[subrefformat=parens,labelformat=parens]{subcaption}
\usepackage{xfrac}

\usepackage[backend=bibtex, style=alphabetic]{biblatex}
\usepackage[noabbrev, capitalise, nameinlink]{cleveref}
\crefname{equation}{}{}
\crefname{assumption}{Assumption}{Assumptions}
\crefformat{equation}{\textup{#2(#1)#3}}

\newtheorem{theorem}{Theorem}[section]
\newtheorem{corollary}[theorem]{Corollary}
\newtheorem{lemma}[theorem]{Lemma}

\theoremstyle{definition}

\theoremstyle{remark}

\numberwithin{theorem}{section}
\numberwithin{equation}{section}
\numberwithin{figure}{section}

\def\TH{\mathcal T_H}

\def\with{\,:\,}

\def\Nb{\mathsf{N}}

\newcommand{\scalar}[2]{({#1},{#2})}

\numberwithin{equation}{section}
\numberwithin{theorem}{section}

\AtBeginDocument{%
	\def\MR#1{}
}

\title[An LOD method for heterogeneous mixed-dimensional problems]{A Localized Orthogonal Decomposition method for heterogeneous mixed-dimensional problems}
\author[M.~Hauck, A.~Målqvist, M.~Mosquera]{Moritz Hauck$^*$, Axel Målqvist$^\dagger$, Malin Mosquera$^\dagger$}
\thanks{${}^*$ Institute for Applied and Numerical Mathematics, Karlsruhe Institute of Technology, Englerstr.~2, 76131 Karlsruhe, Germany, (\upshape{moritz.hauck@kit.edu})}
\thanks{${}^{\dagger}$ Department of Mathematical Sciences, Chalmers University of Technology and University of 	Gothenburg, Göteborg 412 96, Sweden, (\upshape{axel@chalmers.se}, \upshape{malinni@chalmers.se})}

\begin{document}
	
	\begin{abstract}
		We propose a multiscale method for mixed-dimensional elliptic problems with highly heterogeneous coefficients arising, for example, in the modeling of fractured porous media. The method is based on the Localized Orthogonal Decomposition (LOD) framework and constructs locally supported, problem-adapted basis functions on a coarse mesh that does not need to resolve the coefficient oscillations. These basis functions are obtained in parallel by solving localized fine-scale problems. Our a priori error analysis shows that the method achieves optimal convergence with respect to the coarse mesh size, independent of the coefficient regularity, with an exponentially decaying localization error. Numerical experiments validate these theoretical findings and demonstrate the computational viability of the method.
	\end{abstract}
	


\maketitle

\section{Introduction}
The simulation of physical processes in highly heterogeneous media is challenging, especially when the material contains thin structures such as fractures, cracks, or reinforcements. Resolving these features explicitly in large scale simulation can be prohibitively expensive, potentially exceeding the available computing resources. A common and effective alternative is to model these thin structures as lower-dimensional interfaces, embedded within a higher-dimensional bulk domain. Such problems arise in many applications, most prominently in the modeling of flow through fractured porous media. In this context, the fractures are modeled as lower-dimensional manifolds with effective interface conditions; see, e.g.,~\cite{Martin2005,Arrars2019,Boon2018,DelPra2017,Fumagalli2018}. Related problems include modeling vascular networks embedded in biological tissue; see~\cite{Fritz2022}.

From a numerical perspective, solving mixed-dimensional PDEs presents several challenges. Classical finite element methods require meshes that resolve the complex geometry of embedded features, which can become intractable for a large number of interfaces. The situation is further complicated by the presence of strongly heterogeneous and oscillatory coefficients, which must also be resolved by the considered meshes. This situation is typical in porous media flow where the permeability varies rapidly in the bulk regions between interfaces. This paper focuses precisely on such a setting: we study a mixed-dimensional, elliptic diffusion problem with highly heterogeneous coefficients. For such problems, we develop a multiscale method that constructs problem-adapted, locally computable basis functions defined on a coarse mesh. Since these basis functions capture the essential features of the underlying problem, accurate approximations can be obtained even on very coarse meshes that do not resolve the  fine-scale features of the problem.

In general, numerous numerical methods have been developed to solve partial differential equations on surfaces. A foundational contribution to the finite element approximation of the Laplace–Beltrami equation was given in \cite{Dziuk1988}, where the method is formulated on a polyhedral approximation of the surface (see also the comprehensive overview in \cite{Dziuk2013}). Alternative approaches include trace-based methods~\cite{Burman2018,Burman2014,Olshanskii2009}, in which surface functions are represented as traces of functions defined in a higher-dimensional ambient space.
Finite element methods for coupled bulk–interface problems have also been the subject of extensive research; see, for example,~\cite{Elliott2012,Martin2005}. More recently, a fitted finite element method  has been introduced in \cite{Hellman2023}, with particular emphasis on establishing Poincaré-type inequalities for the coupled bulk–interface model.

To address the challenges posed by heterogeneous coefficients in elliptic model problems, a wide range of multiscale methods have been developed. These methods aim to construct low-dimensional problem-adapted approximation spaces that effectively capture fine-scale features. The corresponding basis functions are typically associated with mesh entities of a coarse mesh and computed by solving local fine-scale problems. Prominent approaches include the Heterogeneous Multiscale Method (HMM)\cite{EE03,EE05}, (Generalized) Multiscale Finite Element Methods (GMsFEM) \cite{BabO83,BabCO94}, Multiscale Spectral Generalized Finite Element Methods \cite{BabL11,MaSD22}, the Localized Orthogonal Decomposition (LOD) method~\cite{MalP14,HenP13}, and gamblets \cite{Owh17,OwhS19}.  More recently, refined localization strategies within the LOD framework have been proposed; see, e.g.,~\cite{HauP23,FreHKP24}. For a comprehensive survey of these methods, we refer to~\cite{AltHP21}.
In \cite{Hellman2021}, an initial step toward extending the LOD method to mixed-dimensional models was taken. However, that approach enforces the continuity of the solution across interfaces, making it more restrictive than the setting considered here. Moreover, the methods differ fundamentally: our construction is explicitly tailored to domains with a fractured structure. This allows us to relax the geometric assumptions on the interface and to improve the localization properties of the method.

This paper proposes a multiscale method, in the spirit of the LOD method, for the mixed-dimensional model problem with heterogeneous coefficients. Similar to the classical elliptic setting, the method decomposes the underlying energy space into an infinite-dimensional fine-scale space characterized by the vanishing of integrals with respect to coarse mesh entities, and a finite-dimensional complement space. Crucially, the complement space is defined with respect to the energy inner product, thereby adapting it to the multiscale problem at hand. The canonical basis functions of the complement space are associated with coarse mesh entities and exhibit exponential decay with respect to the coarse mesh. This motivates their localization to patches of elements of the coarse mesh. The precise way the localized basis functions are defined is inspired by the recent work~\cite{HLM25} in the classical elliptic setting. 
We provide a rigorous a priori error analysis of the proposed method in the mixed-dimensional setting, showing that the error consists of two terms. The first term is an optimal-order term in the coarse mesh size. The second term decays exponentially with the number of element layers of the patches used for localization. Thus, increasing the number of element layers logarithmically with the desired accuracy yields an optimally convergent approximation.
Numerical experiments are included to validate and support the theoretical results.

The remainder of the paper is organized as follows. In \cref{sec:modelproblem}, we introduce the mixed-dimensional model problem. A prototypical multiscale method that achieves optimal convergence rates under minimal structural assumptions is presented in \cref{sec:proto}. In \cref{sec:decay}, we establish the exponential decay of the corresponding prototypical basis functions, motivating their localization, which is discussed in \cref{sec:localization}. Based on the resulting localized basis functions, a practical multiscale method is developed in \cref{sec:pracmethod}. In \cref{sec:complexelements}, we extend the methodology to allow for interfaces that vary of the fine computational mesh. Finally, in \cref{sec:numerics}, we present a series of numerical experiments that support our theoretical findings.

\section{Mixed-dimensional model problem}\label{sec:modelproblem}

This section introduces the mixed-dimensional model problem considered in this paper. We begin by specifying the underlying mixed-dimensional geometry.
\subsection{Geometry}
Consider a Lipschitz domain $\Omega \subset \mathbb{R}^d$, $d \in \{2, 3\}$, partitioned into $d+1$ not necessarily connected subdomains of varying dimensionalities as follows:
\begin{equation}
	\Omega = \Omega^0 \cup \cdots \cup \Omega^d,
\end{equation}
where $c$ in $\Omega^c$ denotes the codimension ($0 \leq c \leq d$). Each $\Omega^c$ is decomposed into finitely many disjoint and connected subdomain segments indexed by a set $L_c$, as
\begin{equation}
	\Omega^c = \bigcup_{\ell \in L_c} \Omega^c_\ell,
\end{equation}
where each $\Omega^c_\ell$ is contained in a planar hypersurface of codimension~$c$.

\begin{figure}[h!]
	\centering
	\begin{tikzpicture}[scale=5]
		\draw (0,0) -- (0,1) -- (1,1) -- (1,0) -- (0,0); 
		\draw (0,.2) -- (1,.8); \draw (.4,1) -- (.5,.5); 
		\node at (.8,.2) {$\Omega_1^0$};
		\node at (.7,.8) {$\Omega_2^0$};
		\node at (.2,.7) {$\Omega_3^0$};
		\draw [->] (-.2,.4) to [out=75,in=180] (-.01,.5); \node [below] at (-.2,.4) {$\partial \Omega$};
		\node [circle, fill, black, inner sep=1pt] at (.5,.5) {};
		\draw [<-] (.43,.9) to [out=80,in=200] (.55,1.1); \node [right] at (.55,1.1) {$\Omega^1_1$};
		\draw [<-] (.85,.7) to [out=-30,in=210] (1.1,.77); \node [right] at (1.1,.8) {$\Omega^1_2$};
		\draw [<-] (.25,.33) to [out=-75,in=20] (.19,.22); \node [left] at (.2,.2) {$\Omega^1_3$};
		\draw [<-] (.52,.485) to [out=-25,in=180] (.65,.45); \node [right] at (.65,.45) {$\Omega^2_1$};
	\end{tikzpicture}
	\caption{An example satisfying our assumptions for \( d = 2 \), with subdomain segments of codimensions \( c = 0, 1, 2 \) indicated.
	}\label{fig:domain_rules}
\end{figure}
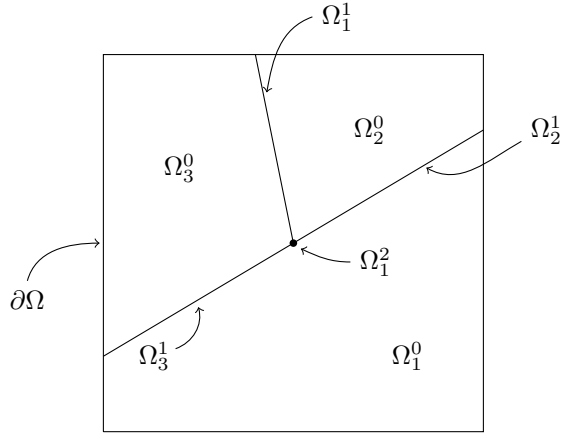

To formulate assumptions on the topology of the problem, we define, for each codimension \( c \), a topological subspace \( X^c \subset \mathbb{R}^d \).  For codimension \( c = 0 \), we set \( X^0 =~\Omega \), and for codimensions \( c > 0 \), we define
\(
X^c = X^{c-1} \setminus \Omega^{c-1}.
\) Standard topological concepts such as openness, closedness, and density are understood with respect to the subspace topology on \( X^c \). The closure, boundary, and interior operators are denoted by \(\overline{\cdot}\), \(\partial \cdot\), and \(\operatorname{int}(\cdot)\), respectively.
We assume that each subdomain segment~\(\Omega^c_\ell\), for \(\ell \in L_c\), is open in \( X^c \), and that their union \(\Omega^c\) is dense in \( X^c \).
Furthermore, we suppose that for any pair \((i, j) \in L_c \times L_{c+1}\), we either have the inclusion \(\Omega^{c+1}_j \subseteq \partial \Omega^c_i\) or it holds  \(\Omega^{c+1}_j \cap \overline{\Omega^c_i} = \emptyset\). 
This condition allows us to define an adjacency relation \( E^c \) between segments of codimensions \( c \) and \( c+1 \) as
\begin{equation}\label{eq:graphEc}
	(i, j) \in E^c \quad \text{if} \quad \Omega^{c+1}_j \subseteq \partial \Omega^c_i,
\end{equation}
and to represent boundary integrals over subdomain segments as sums of subdomain segments of one higher codimension.
In what follows, we focus on subdomains of codimensions 0, 1, and 2.
To simplify the considered setting, we assume that all  subsets \(\Omega^c_\ell\) have Lipschitz boundaries and are polyhedra, polygons, line segments, or points, according to their dimension. This assumption excludes pathological cases such as slit domains and ensures, for example, that an interface cannot end within a bulk region without connecting to another interface.
For an example of a domain satisfying these requirements, see~\cref{fig:domain_rules}.

\subsection{Function spaces}
For a measurable set \( \omega \subset \mathbb{R}^d \), we denote by \( L^2(\omega) \) the Lebesgue space of square-integrable functions on \( \omega \), and by \( H^1(\omega) \) the Sobolev space consisting of functions in \( L^2(\omega) \) whose weak derivatives also belong to \( L^2(\omega) \). The \( L^2(\omega) \)-inner product is denoted by~\( (\cdot, \cdot)_\omega \), with induced norm \( \| \cdot \|_{L^2(\omega)} \). The $H^1(\omega)$-norm is denoted by \( \| \cdot \|_{H^1(\omega)} \).
In the case that the subdomain \( \omega \) is a subdomain segment of codimension \( c > 0 \), we interpret~\( H^1(\omega) \) as \( H^1(\widehat{\omega}) \), where \( \widehat{\omega} \subset \mathbb{R}^{d-c} \) is the image of \( \omega \) under a rigid coordinate transformation. This allows the use of standard tools such as differential operators, the Poincaré inequality, Green's theorem, and trace theorems on  \( \mathbb{R}^{d-c} \). For clarity, we denote the gradient on subdomain segments of codimension $c$ by \( \nabla^c \) with the convention \( \nabla=\nabla^0 \).

We introduce the product space
\[
V \coloneqq V^0 \times V^1,
\]
using the classical notation for product spaces. Here, the space \( V^0 \), defined on the bulk domain of codimension 0, is broken along the interfaces between subdomain segments. Specifically we let
\[
V^0 \coloneqq \prod_{i \in L^0} H^1(\Omega_i^0),
\]
recalling that \(\Omega^0\) is the disjoint union of the subdomains \(\{\Omega_i^0\}_{i \in L^0}\).
The space~$V^1$ is defined on the interfaces of codimension~1. 
Continuity is enforced across shared boundaries, which are formed by subdomain segments of codimension~2, so that
\[
V^1 \coloneqq \Big\{ v \in \prod_{j \in L^1} H^1(\Omega_j^1) \with v_j|_{\Omega_k^2} = v_{j'}|_{\Omega_k^2} \text{ for all } (j, k), (j', k) \in E^1 \Big\},
\]
where the restriction of a function to a subdomain segment is indicated with a subscript index for the respective subdomain.

The subspaces of $V^c$ for $c \in \{0,1\}$, which satisfy homogeneous Dirichlet boundary conditions on $\partial \Omega$ are defined~as
\[
V^c_0 \coloneqq \{ v^c \in V^c \with v^c|_{\partial \Omega} = 0 \},
\]
and their product space is denoted by
\[
V_0 \coloneqq V^0_0 \times V^1_0.
\]

In the following, for a function \(v \in V\), we write \(v = (v^0, v^1)\) and equip the space~\(V\), and naturally the subspace \(V_0\), with the product norm
\begin{equation}
	\label{eq:norm}
	\|v\|_V^2 \coloneqq  \|v^0\|_{H^1(\Omega^0)}^2 + \|v^1\|_{H^1(\Omega^1)}^2.
\end{equation}

\subsection{Problem formulation}

The mixed-dimensional model problem involves coefficients defined on both the bulk domain and the interface domain. On the bulk domain $\Omega^0$, we consider a coefficient $A^0 \in L^\infty(\Omega^0)$. Similarly, on the interface domain $\Omega^1$, we introduce coefficients $A^1 \in L^\infty(\Omega^1)$ and $B^1 \in L^\infty(\Omega^1)$. These coefficients are assumed to satisfy uniform ellipticity and boundedness conditions: there exist positive constants $0 < \underline{\alpha} \le \overline{\alpha} < \infty$ and $0 < \underline{\beta} \le \overline{\beta} < \infty$ such that
\begin{align}
	\label{eq:coeffbound}
	\begin{split}
		\underline{\alpha} \le A^0(x) &\le \overline{\alpha} \quad \text{a.e. in } \Omega^0,\\
		\underline{\alpha} \le A^1(x) \le \overline{\alpha}, \quad \underline{\beta} &\le B^1(x) \le \overline{\beta} \quad \text{a.e. in } \Omega^1.
	\end{split}
\end{align}
The bilinear form \(a\) of the mixed-dimensional model problem is then defined as
\begin{equation}
	\label{eq:defa}
	\begin{aligned}
		a(v, w) & \coloneqq \scalar{A^0 \nabla^0 v^0}{\nabla^0 w^0}_{\Omega^0} + \scalar{A^1 \nabla^1 v^1}{\nabla^1 w^1}_{\Omega^1} \\&
		\qquad + \sum_{(i, j) \in E^0} \scalar{B^1 (v^0_i - v^1_j)}{w^0_i - w^1_j}_{\Omega^1_j}.
	\end{aligned}
\end{equation}

Given a right-hand side functional
\begin{equation}
	\label{eq:rhs}
	F(w) \coloneqq  \scalar{f^0}{w^0}_{\Omega^0} +  \scalar{f^1}{w^1}_{\Omega^1}, 
\end{equation}
with source terms $f^0 \in L^2(\Omega^0)$, $f^1 \in L^2(\Omega^1)$, the mixed-dimensional model problem is formulated as follows: find $u \in V_0$ such that
\begin{equation}
	\label{eq:weak_problem}
	a(u, v) = F(v), \quad \text{for all }v \in V_0.
\end{equation}
The bulk and interface domains are coupled through a Robin-type boundary condition. Continuity of the solution across interface segments is enforced, and a Kirchhoff-type law governs the flux exchange between adjacent interface segments.

If the solution and data are sufficiently smooth, the weak formulation in~\eqref{eq:weak_problem} 
is equivalent to the strong formulation of finding \( u \in V_0 \) such that  
\begin{align*}
	- \nabla \cdot \left( A^0 \nabla u^0 \right) 
	&= f^0,
	&& \text{in } \Omega^0, \\
	- \nabla^1 \cdot \left( A^1 \nabla^1 u^1 \right) 
	- \sum_{(i,j) \in E^0} B^1\!\left(u^0_i - u^1_j\right)
	&= f^1, 
	&& \text{in } \Omega^1,
	\intertext{
		subject to homogeneous Dirichlet conditions on the exterior boundary, continuity of both solution and flux at interface junctions, 
		and the following Robin-type coupling between bulk and interface domains:}
	n_{\Omega^0_i} \cdot \left( A^0 \nabla u^0_i \right) + B^1 u^0_i 
	&= B^1 u^1_j, 
	&& \text{on } \partial \Omega^0_i \cap \Omega^1_j
\end{align*}
for all $i,j$ such that the intersection is non-trivial.
Here, \( n_{\omega} \) denotes the outward unit normal of the subdomain \( \omega \). The solution of this equation can for example model the pressure in a heterogeneous medium in the presence of cracks or faults. 
Further details can be found in~\cite{Hellman2023}.

We close this section by establishing coercivity and boundedness of the bilinear form \(a\), which, since \(V_0\) is a Hilbert space when equipped with the norm \cref{eq:norm}, ensures the well-posedness of problem~\cref{eq:weak_problem} by the Riesz representation theorem.
\begin{lemma}[Properties of $a$]
\label{lem:coercivity}
There exist constants $0<c_a\leq C_a < \infty$, such that
\begin{align*}
	a(v,v) \geq c_a  \|v\|_V^2,\quad v \in V_0,\qquad 
	a(v,w) \leq C_a  \|v\|_V\|w\|_V,\quad v,w \in V_0,
\end{align*}
where, up to constants solely dependent on the geometry of~$\{\Omega^c_\ell\}_{c, \ell}$,  the constants scale as
$c_a \sim (\underline \alpha^{-1} + \underline \beta^{-1})^{-1}$ and $C_a \sim (\overline \alpha + \overline \beta)$.
\end{lemma}

A direct consequence of this result is that \( \| \cdot \|_a^2 \coloneqq a(\cdot, \cdot) \) defines a norm on \( V_0 \), which will henceforth be referred to as the energy norm.

\begin{proof}
The coercivity proof employs an iterative procedure, starting by bounding norms on subdomain segments far from the domain boundary and progressively extending these bounds to neighboring segments until reaching those where a Friedrichs-type inequality applies; see \cite[Lem.~1]{Hellman2023} for details. The boundedness can be proved using a standard trace inequality; see \cite[Lem.~2]{Hellman2023}.
\end{proof}

\section{Prototypical method}\label{sec:proto}
In this section, we present a prototypical multiscale method that achieves optimal-order approximations without any  pre-asymptotic effects, under minimal structural assumptions on the coefficients.  
The LOD method is typically defined with respect to a coarse mesh of the domain. In the present mixed-dimensional setting, coarse meshes are required for both the bulk domain and the interface domain, denoted by \( \mathcal{T}_H^0 \) and \( \mathcal{T}_H^1 \), respectively. They are henceforth referred to as the \emph{bulk mesh} and the \emph{interface mesh}, respectively.  Note that the considered coarse meshes are not required to resolve the fine-scale variations of the coefficients.
In general, there is considerable flexibility in how these meshes are chosen. Even complex-shaped elements representing fractures can be accommodated, provided the elements are non-degenerate and a local element-wise  Poincaré inequality holds.  

Subsequently, we restrict our attention to a significantly simplified setting; see \cref{sec:complexelements} for an extension of the methodology to more general meshes.
Specifically, we consider a bulk mesh that is a quasi-uniform partition of \( \Omega \) into finitely many closed, shape-regular simplicial or quadrilateral/hexahedral elements, with diameter at most \( H > 0 \). We assume that the bulk mesh is such that each bulk subdomain segment~\(\Omega_i^0\) can be expressed as a union of bulk mesh elements. Similarly, the interface mesh is assumed to be such that each interface subdomain segment~\(\Omega_j^1\) can be expressed as a union of interface mesh elements. Moreover, we assume the compatibility condition that any face of a bulk mesh element lying on an interface~\(\Omega_j^1\) coincides exactly with an element of the interface mesh. Consequently, the mesh size \(H\) is also the maximal element diameter of the interface mesh.

Following the presentation of the LOD in~\cite{AltHP21}, the definition of the method is based on so-called \emph{quantities of interest} (QOI). Their purpose is to encode which information about the solution \( u \) to problem~\cref{eq:weak_problem} the prototypical multiscale method should preserve exactly.
Each QOI is associated with an element of either the bulk or the interface mesh and is defined as the corresponding element-wise average:
\begin{equation}
\label{defQOI}
q_{T^c} (v) \coloneqq  \frac{1}{|T^c|_{d-c}}\int_{T^c} v^c \, \mathrm{d}x,\quad T^c \in \mathcal T_H^c,\; c \in \{0,1\},
\end{equation}
where \(|\cdot|_{n}\) for \(n \in \{d-1,d\}\) denotes the \(n\)-dimensional volume. Based on these QOI, we can now define a so-called \emph{fine-scale space} as
\begin{equation}
\label{eq:defW}
W \coloneqq \left\{v \in V_0 \with q_{T^c}(v) = 0 \text{ for all }T^c \in \mathcal T_H^c,\; c \in \{0,1\}\right\}.
\end{equation}
Intuitively, this space consists of functions that average out on length scales smaller than or equal to \( H \). It represents the fine-scale information that cannot be captured on the coarse scale.

The approximation space of the prototypical multiscale method is then defined as the (finite-dimensional) orthogonal complement of the fine-scale space \(W\) in \(V_0\), with respect to the energy inner product \(a\), i.e.,
\begin{equation}
\label{eq:Zms}
\tilde V_H \coloneqq \left\{ v \in V_0 \with  a (v, w) = 0 \text{ for all }w \in W \right\}.
\end{equation}
Since \(\tilde V_H\) is constructed as the orthogonal complement of \(W\) with respect to \(a\), it inherently encodes problem-specific information. This enables reliable approximations even at scales where the coefficients are not resolved. The use of tildes in the notation serves to emphasize that the space is adapted to the problem at hand. A basis for the space \( \tilde V_H \) can be constructed as outlined in the following lemma.

\begin{lemma}[Prototypical basis]\label{le:protbasis}
A basis of the space $\tilde V_H$ is given by
\begin{equation}
	\label{eq:LODprotbasis}
	\big\{\tilde \varphi_{T^c}\with T^c \in \mathcal T_H^c,\, c \in \{0,1\}\big\}
\end{equation}
with $\tilde \varphi_{T^c}$ defined for all $T^c \in \mathcal T_H^c$, $c \in \{0,1\}$ as the unique solutions to the saddle point problem: find $(\tilde \varphi_{T^c},\lambda) \in V_0\times \mathbb R^N$ with $N \coloneqq \#\mathcal T_H^0+\#\mathcal T_H^1$ such that
\begin{subequations}
	\label{pbphiE} 
	\begin{align}
		&\qquad \qquad a (\tilde \varphi_{T^c}, v)& +&  &b(v,\lambda) & &=\quad  &0, &&\text{ for all }v \in V_0, &&\qquad \qquad \label{eq:LODbasis1}\\
		&\qquad\qquad b(\tilde \varphi_{T^c},\mu)                   &   &         &    & &=\quad  &\mu_{T^c}, &&\text{ for all }\mu \in \mathbb R^N,&&\qquad \qquad  \label{eq:LODbasis3}
	\end{align}
\end{subequations}
where the bilinear form $b \colon V \times \mathbb R^{N}\to \mathbb R$ is defined as
\begin{equation*}
	b(v,\mu) \coloneqq \sum_{c \in \{0,1\}}\sum_{T^c \in \mathcal T_H^c} \mu_{T^c}  q_{T^c}(v),
\end{equation*}
and \(\mu_{T^c}\) denotes the component of the vector \(\mu\) associated with the element \(T^c\).
\end{lemma}

Before proving this lemma, we introduce the notion of bubble functions, which will be used not only in this proof but also throughout the manuscript. These are locally supported functions associated with either elements of the bulk or interface meshes and are denoted by \( b_{T^c} \) for \( T^c \in \mathcal{T}_H^c \), \( c \in \{0,1\} \). They are defined as \( b_{T^0} \coloneqq (v_{T^0}, 0) \) for \( T^0 \in \mathcal{T}_H^0 \) and \( b_{T^1} \coloneqq (0, v_{T^1}) \) for \( T^1 \in \mathcal{T}_H^1 \), where \( v_{T^c} \in H^1_0(T^c) \) is chosen such that the resulting bubble function \( b_{T^c} \) satisfies the Kronecker-delta property:
\begin{equation}\label{eq:kdbubble}
q_{K^{c'}}(b_{T^c}) = \delta_{K^{c'} T^c},\quad 
K^{c'} \in \mathcal{T}_H^{c'},\;  c' \in \{0,1\},
\end{equation}
and such that the following stability estimate holds:
\begin{align}\label{eq:stabbubble}
\begin{split}
	\|\nabla^c v_{T^c}\|_{L^2(T^c)} &\lesssim H^{-1} \|v_{T^c}\|_{L^2(T^c)} \lesssim H^{(d-c)/2 - 1},\quad c \in \{0,1\}.
\end{split}
\end{align}
The construction of the functions $v_{T^c}$ is classical; they can be obtained as suitably scaled products of hat functions. We emphasize that bubble functions only serve as a theoretical tool and are never used explicitly in practical implementations.

\begin{proof}[Proof of \cref{le:protbasis}]
We begin by establishing well-posedness of problem~\cref{pbphiE} by verifying the inf--sup stability of the bilinear form~$b$. Specifically, we show that there exists a constant $C > 0$ such that
\begin{equation}\label{eq:infsupdoublespp}
	\adjustlimits \inf_{\mu \in \mathbb{R}^N} \sup_{v \in V} \frac{|b(v,\mu)|}{\|v\|_V \,|\mu|} \geq C,
\end{equation}
where $|\cdot|$ denotes the Euclidean norm on~$\mathbb{R}^N$.

Given any $\mu \in \mathbb{R}^N$, we construct a function $v \in V_0$ componentwise as  
\begin{align*}
	v \coloneqq  \sum_{c \in \{0,1\}}\sum_{T^c \in \mathcal T_H^c} \mu_{T^c}\, b_{T^c}
\end{align*}
using the bubble functions introduced above.  
Noting that, by construction, we have $|b(v, \mu)| = |\mu|^2$, and employing the locality of the bubble functions together with a local Friedrichs inequality and the stability estimates from~\cref{eq:stabbubble}, we obtain
\begin{align*}
	\|v\|_V^2 &\lesssim \sum_{c \in \{0,1\}}\sum_{T^c \in \mathcal T_H^c} |\mu_{T^c}|^2 \|\nabla^c v_{T^c}\|_{L^2(T^c)}^2 
	\lesssim H^{d-2} |\mu|^2.
\end{align*}
Combining these estimates establishes the inf–sup condition~\cref{eq:infsupdoublespp}.	

Next, we prove that~\cref{eq:LODprotbasis} defines a basis for the space~$\tilde V_H$. To this end, we first observe that condition~\cref{eq:LODbasis1} implies that $\tilde \varphi_{T^c} \in \tilde V_H$ for all $T^c \in \mathcal T_H^c$, $c \in \{0,1\}$. 
To establish the basis property, let $v \in \tilde V_H$ be arbitrary, and define
\[
w \coloneqq v - \sum_{c \in \{0,1\}}\sum_{T^c \in \mathcal T_H^c}q_{T^c}(v)\tilde \varphi_{T^c}.
\]
By construction, it holds that $w \in W$, and since $w$ is also $a$-orthogonal to $W$, it follows that $w = 0$. Thus, any $v \in \tilde V_H$ can be expressed uniquely as a linear combination of functions in \cref{eq:LODprotbasis}, which proves that they form a basis of~$\tilde V_H$. This completes the proof.
\end{proof}
Having introduced the prototypical basis functions, we can define a projection operator $\mathcal R \colon V \to \tilde V_H$ that preserves the QOI from~\cref{defQOI} by
\begin{equation}
\label{eq:defR}
\mathcal R v \coloneqq \sum_{c \in \{0,1\}}\sum_{T^c \in \mathcal T_H^c}  
q_{T^c}(v)\tilde \varphi_{T^c},\quad v \in V_0.
\end{equation}
This operator coincides with the $a$-orthogonal projection onto $\tilde V_H$, since for any $v \in~V_0$ and $w \in \tilde V_H$ it holds that $a(v-\mathcal R v,w)=0$, noting that $v-\mathcal R v\in W$. This orthogonality relation also implies the continuity of  $\mathcal R$, i.e.,
\begin{equation}
\label{eq:Rcont}
\|\mathcal R v\|_V \leq \sqrt{C_a/c_a} \|v\|_V,\quad v \in V_0,
\end{equation}
where we have used \cref{lem:coercivity}.

The desired prototypical multiscale method then seeks $\tilde u_H \in \tilde V_H$ such that
\begin{align}\label{eq:protmethod}
a(\tilde u_H, \tilde v_H) = F(\tilde v_H) \quad \text{for all } \tilde v_H \in \tilde V_H,
\end{align}
where $F$ is defined in~\cref{eq:rhs}. The following theorem proves optimal-order convergence of this method under minimal structural assumptions.

\begin{theorem}[Prototypical method]\label{thm:convergenceprot}
The prototypical multiscale method~\cref{eq:protmethod} is well-posed, and its solution is given by $\tilde u_H = \mathcal{R}u$. Consider any right-hand side~$F$ of the form \cref{eq:rhs} with components satisfying $f^0 \in H^s(\Omega^0)$ and $f^1 \in H^s(\Omega^1)$, where $s \in \{0,1\}$, and define a corresponding norm by
\begin{equation}\label{eq:fnorm}
	|f|_s^2 \coloneqq  |f^0|_{H^s(\Omega^0)}^2 +  |f^1|_{H^s(\Omega^1)}^2,
\end{equation}
where $|\cdot|_{H^s}$ denotes the $H^s$-seminorm. Then, the following error estimate holds:
\begin{align}
	\|u - \tilde u_H\|_V &\lesssim H^{1+s} |f|_s. \label{eq:errestH1}
\end{align}
\end{theorem}
\begin{proof}
First, observe that the space $\tilde V_H$ is a closed subspace of $V_0$. 
Hence, by the Riesz representation theorem, problem~\cref{eq:protmethod} is well-posed. 
Comparing \eqref{eq:Zms} with the weak formulation \eqref{eq:weak_problem}, we obtain that
\begin{equation}
	\label{eq:orthogonalityrelation}
	a(u - \tilde u_H, \tilde v_H) = 0,\quad \tilde v_H \in \tilde V_H,
\end{equation}
which shows that $\tilde u_H = \mathcal R u$ with the $a$-orthogonal projection $\mathcal R$ introduced above.

To derive the error estimate, we observe that $e \coloneqq u - \tilde u_H$ lies in the space~$W$. Using the coercivity result from~\cref{lem:coercivity}, the orthogonality relation \cref{eq:orthogonalityrelation}, and testing the weak formulation~\cref{eq:weak_problem} with $e$, we obtain
\begin{equation}
	\label{eq:errestprot}
	c_a \|e\|_V^2 \leq a(e, e) = a(u, e) = F(e).
\end{equation}
Let $\Pi_H^0 \colon L^2(\Omega^0) \to \mathbb{P}^0(\mathcal{T}_H^0)$ and $\Pi_H^1 \colon L^2(\Omega^1) \to \mathbb{P}^0(\mathcal{T}_H^1)$ denote the $L^2$-orthogonal projections onto piecewise constant functions on the coarse grid. Exploiting their orthogonality properties yields the following representation for $F(e)$:
\begin{equation*}
	F(e) =  (f^0 - \Pi_H^0 f^0,\, e^0 - \Pi_H^0 e^0)_{\Omega^0} + (f^1 - \Pi_H^1 f^1,\, e^1 - \Pi_H^1 e^1)_{\Omega^1},
\end{equation*}
where the $L^2$-projections of $f^0$ and $f^1$ may be omitted without affecting the identity.
Estimate~\cref{eq:errestH1} then follows by applying a local Poincaré-type inequality.
\end{proof}

\section{Exponential decay of prototypical basis}\label{sec:decay}

We emphasize that the prototypical LOD basis functions defined in~\cref{eq:LODprotbasis} are globally supported, and thus their computation requires solving global problems. However, this would be computationally infeasible in practice.  In this section, we demonstrate that these basis functions exhibit an exponential decay, justifying their approximation by localized counterparts defined on oversampling domains, as will be discussed  in~\cref{sec:localization}. A practical multiscale method based on these localized basis functions will then be introduced in~\cref{sec:pracmethod}.

For the localization procedure, we adopt the approach proposed in the recent generalized framework for high-order LOD methods; see \cite{HLM25}. This framework accommodates the choice of QOI on the bulk domain made here and enables a stable basis construction, ensuring that the error of the localized method does not increase as the coarse mesh is refined, when the oversampling parameter (the number of element layers composing the oversampling domains) remains fixed. Such negative effects have been observed, for example, in \cite{MalP14,Maier2021}. Notably, this stability is achieved without relying on intricate bubble function constructions in the practical implementation of the method, cf. \cite{Hauck2022,Dong2023}.
The localization is based on a quasi-interpolation operator $\mathcal{I}_H \colon V \to V_H$, where $V_H \coloneqq V_H^0 \times V_H^1$, with $V_H^0$ and $V_H^1$ denoting the conforming piecewise linear finite element spaces defined with respect to the bulk and interface meshes, respectively. In particular, we have $V_H \subset V$. The operator is defined component-wise as $\mathcal{I}_H v = (\mathcal{I}_H^0 v, \mathcal{I}_H^1 v)$, where each component $\mathcal{I}_H^cv$ is defined at interior nodes $z^c$ of the mesh $\mathcal{T}_H^c$ for codimension $c \in \{0,1\}$ as:
\begin{align}
\label{eq:IH}
(\mathcal{I}_H^c v)(z^c) &\coloneqq \frac{1}{\# \mathcal{T}^c(z^c)} \sum_{T^c \in \mathcal{T}^c(z^c)} q_{T^c}^c(v^c) \Lambda_{z^c},
\end{align}
where $\mathcal{T}^c(z^c)$ denotes the set of elements in $\mathcal{T}_H^c$ sharing the node $z^c$, and $\Lambda_{z^c}$ is the classical hat function associated with node $z^c$. For boundary nodes, the corresponding nodal values of $\mathcal{I}_H^c v$ are set to zero to conform to the prescribed Dirichlet boundary conditions.
For $c \in \{0,1\}$, this operator satisfies the following local approximation and stability estimates:
\begin{align}
\label{eq:propIH}
H^{-1}\|v - \mathcal{I}_H^c v\|_{L^2(T^c)} + \|\nabla^c \mathcal{I}_H^c v\|_{L^2(T^c)} \lesssim \|\nabla^c v\|_{L^2(\mathsf{N}^c(T^c))},\; v \in V^c,\; T^c \in \mathcal T_H^c
\end{align}
with $\mathsf{N}^c(S^c)$ denoting the first-order element patch of a union of elements $S^c$ in $\mathcal{T}_H^c$, defined as
\begin{equation}
\label{eq:patchesc}
\mathsf{N}^c(S^c) \coloneqq \bigcup \left\{ T^c \in \mathcal{T}_H^c \colon \overline{S^c} \cap \overline{T^c} \neq \emptyset \right\}.
\end{equation}
Note that constructions such as \cref{eq:IH} and the properties in \cref{eq:propIH} are classical; see, e.g.,~\cite{ErnG17}. The specific choice of quasi-interpolation operator is not essential. Apart from the stability estimate~\eqref{eq:propIH}, we only require that \(\mathcal{I}_H v = 0\) holds for all \(v \in W\).

Following the localization strategy in \cite{HLM25}, the $a$-orthogonal projection operator $\mathcal{R}$ defined in \cref{eq:defR} is decomposed as:
\begin{equation}
\mathcal R = \mathcal I_H - \mathcal K,
\end{equation}
where $\mathcal{K} \colon V_0 \to V_0$ maps any $v \in V_0$ to the unique solution $(\mathcal{K}v, \lambda) \in V_0 \times \mathbb{R}^N$ of
\begin{subequations}
\label{eq:defK} 
\begin{align}
	&\quad \quad a (\mathcal Kv, w)& +&  &b(w,\lambda) & &=\quad  &a(\mathcal I_H v,w), &&\text{ for all }w \in V_0, &&\quad \label{eq:defK1}\\
	&\quad\quad b(\mathcal Kv,\mu)                   &   &         &    & &=\quad  &-b(v-\mathcal I_H v,\mu), &&\text{ for all }\mu \in \mathbb R^N.&& \quad\label{eq:defK2}
\end{align}
\end{subequations}
Note that \(\lambda\) is a generic notation for a Lagrange multiplier; the specific multiplier being referred to will be clear from the context.

The operator $\mathcal{K}$ can be decomposed as
\begin{equation*}
\mathcal{K} = \sum_{T^0 \in \mathcal{T}_H^0} \mathcal{K}_{T^0},
\end{equation*}
where each local operator $\mathcal{K}_{T^0} \colon V_0 \to V_0$ maps any $v \in V_0$ to the unique solution $(\mathcal{K}_{T^0} v, \lambda) \in V_0 \times \mathbb{R}^N$ of
\begin{subequations}
\label{eq:defKT} 
\begin{align}
	& a (\mathcal K_{T^0}v, w)& +&  &b(w,\lambda) & &=\quad  &a_{T^0}(\mathcal I_H v,w), &&\text{ for all }w \in V_0, &&\quad \label{eq:defKT1}\\
	& b(\mathcal K_{T^0}v,\mu)                   &   &         &    & &=\quad  &-b_{T^0}(v-\mathcal I_H v,\mu), &&\text{ for all }\mu \in \mathbb R^N,&& \quad\label{eq:defKT2}
\end{align}
\end{subequations}
where \( a_{T^0} \) and \( b_{T^0} \) denote, for any \( T^0 \in \mathcal{T}_H^0 \), suitable restrictions of the bilinear forms \( a \) and \( b \), satisfying the properties
\[
\sum_{T^0 \in \mathcal{T}_H^0} a_{T^0} = a \quad \text{and} \quad \sum_{T^0 \in \mathcal{T}_H^0} b_{T^0} = b.
\]
Denoting by $
n(T^1) \coloneqq \#\{ K^0 \in \mathcal{T}_H^0 \with \partial K^0 \supset T^1\}$ the number of bulk elements sharing the interface element \( T^1 \in \mathcal{T}_H^1 \), such restrictions \( a_{T^0} \) and \( b_{T^0} \) can be defined~as:
\begin{align}
\begin{split}\label{eq:defAT0}
	a_{T^0}(v, w) &\coloneqq  \scalar{A^0 \nabla^0 v^0}{\nabla^0 w^0}_{T^0} + \sum_{T^1 \subset \partial T^0} \frac{1}{n(T^1)} \scalar{A^1 \nabla^1 v^1}{\nabla^1 w^1}_{T^1} \\
	&\qquad+ \sum_{T^1 \subset \partial T^0} \scalar{B^1 (v^0|_{T^0} - v^1)}{w^0|_{T^0} - w^1}_{T^1},
\end{split}
\\
b_{T^0}(v,\mu) &\coloneqq \mu_{T^0}q_{T^0}(v) + \sum_{T^1 \subset \partial T^0} \frac{1}{n(T^1)}\mu_{T^1}q_{T^1}(v).
\end{align}

In what follows, we establish the exponential decay of the operators $\mathcal{K}_{T^0}$ away from the element $T^0$. To quantify the decay, we introduce a notion of patches with respect to the coarse mesh $\mathcal{T}_H$. For an oversampling parameter $\ell \in \mathbb{N}$ and a union of elements $S^c \subset \mathcal{T}_H^c$ with $c \in \{0,1\}$, the $\ell$-th order patch is defined recursively by
\[
\mathsf{N}_\ell^c(S^c) \coloneqq \mathsf{N}_1^c(\mathsf{N}_{\ell-1}^c(S^c)), \quad \ell \geq 2,
\]
with $\mathsf{N}_1^c(S^c) \coloneqq \mathsf{N}^c(S^c)$, as defined in \cref{eq:patchesc}.
We also introduce a notation that jointly represents the bulk and interface components, namely $\mathsf{N}_\ell(S) \coloneqq (\mathsf{N}_\ell^0(S^0), \mathsf{N}_\ell^1(S^1))$ for the $\ell$-th order patches, and $\Omega \setminus \mathsf{N}_\ell(S) \coloneqq (\Omega^0 \setminus \mathsf{N}_\ell^0(S^0), \Omega^1 \setminus \mathsf{N}_\ell^1(S^1))$ for their complement. Additionally, we abbreviate $\mathsf{N}_\ell(T^0) \coloneqq \mathsf{N}_\ell(T^0, \partial T^0)$.

To formulate the desired exponential decay result, we require localized versions of the $V$-norm defined in \cref{eq:norm}. Instead of naively localizing it by merely restricting the $H^1$-norms on the bulk and interface domains to subdomains, we define, for tuples $S = (S^0, S^1)$ as above, the localized norm as follows:
\begin{align}
\label{eq:restanorm}
\begin{split}
	\|v\|_{V(S)}^2&\coloneqq \| \nabla^0v^0\|_{L^2(S^0)}^2+\| \nabla^1v^1\|_{L^2(S^1)}^2 
	+\hspace{-1.25ex} \sum_{\substack{\{T^1\subset \partial T^0\with \\T^0 \subset S^0,\, T^1\subset S^1\}}}\hspace{-1.25ex} \|v^0|_{T^0} - v^1\|_{L^2(T^1)}^2.
\end{split}
\end{align}
This definition will prove more convenient for the subsequent proof. Note that the unrestricted version of this norm (i.e., the above definition with \(S = (\Omega^0, \Omega^1)\)), which we denote by \(\|\cdot\|_{V(\Omega)}\), is equivalent to the energy norm, which, by \cref{lem:coercivity}, is in turn equivalent to the \(V\)-norm defined in \cref{eq:norm} on the space \(V_0\).

\begin{theorem}[Exponential decay]\label{thm:dec}
There exists a constant $c>0$ independent of $H$, $\ell$, and $T^0$ such that for all $T^0 \in \TH$ and $\ell \in \mathbb N$, it holds that
\begin{equation}
	\label{eq:dec}
	\|\mathcal K_{T^0} v\|_{V(\Omega\setminus \Nb_\ell(T^0))} \lesssim \exp(-c \ell)\|\mathcal K_{T^0} v\|_{V(\Omega)},\quad v \in V_0.
\end{equation}
\end{theorem}
\begin{proof}
For a given \( v \in V_0 \), define \( \varphi \coloneqq \mathcal{K}_{T^0} v \). Introduce a cut-off function \( \eta \in V_H^0 \), defined on the bulk domain, and uniquely characterized by the following properties:
\begin{equation}\label{eq:eta}
	\begin{cases}
		\eta(x) =0, & x \in \mathsf{N}_{\ell-1}^0(T^0), \\
		\eta(x) =1, & x \in \Omega^0 \setminus \mathsf{N}_{\ell}^0(T^0), \\
		0 \leq \eta(x) \leq 1, & x \in \mathsf{N}_{\ell}^0(T^0) \setminus \mathsf{N}_{\ell -1}^0(T^0).
	\end{cases}
\end{equation}
By construction the cut-off function satisfies \( \| \nabla^0\eta \|_{L^\infty} \lesssim H^{-1} \). 
Let the ring of elements containing the support of the gradient of $\eta$ be denoted by
\begin{equation}
	\label{eq:ring}
	R_{T^0} =(R_{T^0}^0,R_{T^0}^1)\coloneqq \mathsf{N}_{\ell}(T^0) \setminus \mathsf{N}_{\ell-1}(T^0).
\end{equation}
In what follows, we assume that \(\ell \geq 2\), which implies that \(\overline{T^0} \cap \overline{R_{T^0}^0} = \emptyset\).
For all elements \( T^0 \in \mathcal{T}_H^0 \) and \( T^1 \in \mathcal{T}_H^1 \) contained in the support of \( \eta \), it thus holds that
\begin{equation}\label{eq:avgvarphi}
	\int_{T^0} \varphi^0\,\mathrm{d}x = \int_{T^1} \varphi^1\,\mathrm{d}\sigma = 0.
\end{equation}
To prove this result, we use the identity \( b(\varphi, \mu) = -b_{T^0}(v - \mathcal{I}_H v, \mu) \), which follows from \cref{eq:defKT2}, together with the locality of the restricted bilinear form \( b_{T^0} \).

Using definitions~\cref{eq:defa,eq:restanorm}, the properties of \( \eta \) from~\cref{eq:eta}, integration by parts, the coefficient bounds from~\cref{eq:coeffbound}, and~\cref{eq:avgvarphi}, we obtain that
\begin{align}
	\|\varphi&\|^2_{V(\Omega\setminus \mathsf N_\ell(T^0))} \notag\\
	&\leq
	(\eta A^0  \nabla^0 \varphi^0,\nabla^0 \varphi^0)_{L^2(\Omega^0)} 
	+ (\eta A^1 \nabla^1 \varphi^1,\nabla^1 \varphi^1)_{L^2(\Omega^1)} \notag \\
	&\quad + \sum_{(i, j) \in E_0}(\eta B^1(\varphi^0_i - \varphi^1_j),\varphi^0_i - \varphi^1_j)_{L^2(\Omega^1_j)} \notag \\
	&= 
	a(\varphi,\eta\varphi)
	- (A^0 \nabla^0 \varphi^0, (\nabla^0 \eta) \varphi^0)_{L^2(\Omega^0)} 
	- (A^1 \nabla^1 \varphi^1, (\nabla^1 \eta) \varphi^1)_{L^2(\Omega^1)} \notag \\
	&\lesssim |b(\lambda ,\eta \varphi)| 
	+ \| \nabla^0 \varphi^0\|^2_{L^2(R_{T^0}^0)} 
	+ \|\nabla^1 \varphi^1\|^2_{L^2(R_{T^0}^1)}\label{eq:xi1}.
\end{align}
In the last step, we have tested \cref{eq:defKT1} with \( \eta \varphi \) and used that \( a_{T^0}(\mathcal{I}_H v, \eta \varphi) = 0 \).

We estimate the first term on the right-hand side of \cref{eq:xi1} using the properties of \( \eta \) from \cref{eq:eta}, in combination with \cref{eq:avgvarphi}, yielding
\begin{align}\label{eq:cstart}
	\begin{split}
		|b(\lambda, \eta \varphi)|&=\bigg|\sum_{c \in \{0,1\}}\sum_{K^c \subset R_{T^0}^c} \lambda_{K^c} q_{K^c}(\eta \varphi)\bigg|\\
		&\lesssim \sum_{c \in \{0,1\}}\Big( H^{-d+c+2}\sum_{K^c\subset R_{T^0}^c}  |\lambda_{K^c}|^2\Big)^{1/2}  \|\nabla^c\varphi^c\|_{L^2(R_{T^0}^c)}.
	\end{split}
\end{align}
To derive bounds on the entries of the vector \( \lambda \in \mathbb{R}^N \), we employ element and face bubble functions, as introduced for the proof of \cref{le:protbasis}. By testing \cref{eq:defKT1} with \(b_{K^0} =  (v_{K^0}, 0) \) and using \cref{eq:stabbubble}, we obtain, for any \( K^0 \subset R_{T^0}^0 \), that
\begin{align*}\label{eq:lambdastart}
	|\lambda_{K^0}|&=|b(\lambda, b_{K^0})|=|a(\varphi,b_{K^0})| =|(A^0\nabla^0\varphi^0,\nabla^0v_{K^0})_{L^2(\Omega^0)}|\\
	&\lesssim H^{d/2-1}\|\nabla^0\varphi^0\|_{L^2(K^0)}.
\end{align*}
The second equality follows from the fact that \( a_{T^0}(\mathcal{I}_H v, b_{K^0}) = 0 \).
Summing over all \( K^0 \subset R_{T^0}^0 \) then yields
\begin{equation}\label{eq:sumlamK0}
	H^{-d+2}	\sum_{K^0 \subset R_{T^0}^0} | \lambda_{K^0}|^2 \lesssim \|\nabla^0\varphi^0\|^2_{L^2(R_{T^0}^0)}.
\end{equation}
To derive a similar bound for the entries \( \lambda_{K^1} \) associated with \( K^1 \subset R_{T^0}^1 \), we employ a trace inequality for \( \varphi^0 \) of the form
\begin{equation}
	\label{eq:trace}
	\|\varphi^0\|_{L^2(\partial T^0)} \lesssim H^{-\frac{1}{2}} \|\varphi^0\|_{L^2(T^0)} + H^{\frac{1}{2}} \|\nabla^0 \varphi^0\|_{L^2(T^0)} \lesssim H^{\frac{1}{2}} \|\nabla^0 \varphi^0\|_{L^2(T^0)},
\end{equation}
where the last estimate follows from the Poincaré inequality, which is applicable due to \cref{eq:avgvarphi}. Testing \cref{eq:defKT1} with \(b_{K^1}= (0, v_{K^1}) \), and using that \( a_{T^0}(\mathcal{I}_H v, b_{K^1}) = 0 \) together with stability estimate \cref{eq:stabbubble}, we obtain, for any \( K^1 \subset R_{T^0}^1 \), that
\begin{align*}
	|\lambda_{K^1}|&=|b(\lambda, b_{K^1})|=|a(\varphi,b_{K^1})| \\
	&=|(A^1\nabla^1 \varphi^1,\nabla^1 v_{K^1})_{L^2(K^1)}+\sum_{K^0 \subset \omega_{K^1}}(B(\varphi^0|_{K^0}-\varphi^1),-v_{K^1})_{L^2(K^1)}|\\
	&\lesssim  H^{d/2}\|\nabla^0\varphi^0\|_{L^2(\omega_{K^1})} + H^{(d-1)/2-1}\|\nabla^1 \varphi^1\|_{L^2(K^1)},
\end{align*}
where \( \omega_{K^1} \) denotes the union of the two bulk elements sharing the interface \( K^1 \). Summing over all $K^1 \subset R_{T^0}^1$, we obtain that
\begin{equation}\label{eq:sumlamK1}
	H^{-d+3} \sum_{K^1\subset  R_{T^0}^1} |\lambda_{K^1}|^2 \lesssim \|\nabla^0\varphi^0\|^2_{L^2(\mathsf N^0(R_{T^0}^0))}+\|\nabla^1 \varphi^1\|^2_{L^2(R_{T^0}^1)}.
\end{equation}
Combining \cref{eq:sumlamK0,eq:sumlamK1}, we obtain that
\begin{equation*}
	|b(\lambda,\eta \varphi)|\lesssim \|\nabla^0\varphi^0\|^2_{L^2(\mathsf N ^0(R_{T^0}^0))}+\|\nabla^1 \varphi^1\|^2_{L^2(R_{T^0}^1)}.
\end{equation*}

We conclude that there exists a constant \( C > 0 \), independent of \( H \), \( \ell \), and \( T^0 \), such that the following estimate holds for all $\ell \geq 2$:
\[
\|\varphi\|^2_{V( \Omega \setminus \mathsf{N}_\ell(T^0))} \leq C \|\varphi\|^2_{V( \mathsf{N}_{\ell+1}(T^0) \setminus \mathsf{N}_{\ell-2}(T^0))}.
\]
By the index shift \(\ell \leftarrow \ell - 1\), we obtain, for all \(\ell \geq 3\), that
\begin{align*}
	\|\varphi\|^2_{V( \Omega \setminus \mathsf{N}_\ell(T^0))} 
	&\leq \|\varphi\|^2_{V( \Omega \setminus \mathsf{N}_{\ell-1}(T^0))} \leq C \|\varphi\|^2_{V( \mathsf{N}_\ell(T^0) \setminus \mathsf{N}_{\ell-3}(T^0))} \\
	&= C \left( \|\varphi\|^2_{V( \Omega \setminus \mathsf{N}_{\ell-3}(T^0))} - \|\varphi\|^2_{V( \Omega \setminus \mathsf{N}_\ell(T^0))} \right).
\end{align*}
Rearranging terms and defining \( \gamma \coloneqq \frac{C}{1 + C} < 1 \), we obtain the recursive inequality
\[
\|\varphi\|^2_{V( \Omega \setminus \mathsf{N}_\ell(T^0))} \leq \gamma \|\varphi\|^2_{V( \Omega \setminus \mathsf{N}_{\ell - 3}(T^0))} 
\leq \gamma^{\lfloor \ell / 3 \rfloor} \|\varphi\|^2_{V(\Omega)},
\]
and, noting that there exist constants \( C, c > 0 \) such that $
\gamma^{\lfloor \ell/3 \rfloor / 2} \leq C \exp(-c \ell)$, the assertion can be concluded for $\ell \geq 3$. For \( \ell \in \{1, 2\} \), the result follows from a scaling argument (with a possible different hidden constant).
\end{proof}

\section{Localization}
\label{sec:localization}
The exponential decay properties established in the previous section motivate the localized computation of the prototypical LOD basis functions. The localization strategy proposed here, based on the approach presented in \cite{HLM25}, overcomes a key limitation of the naive localization method, which restricts problem~\eqref{pbphiE} to the \(\ell\)-th order patches. For this naive approach, the localization error deteriorates as the mesh is refined if the parameter \(\ell\) remains fixed; see~\cite{MalP14,Maier2021}.

To formalize the localization strategy, we introduce the localized spaces
\begin{align}
V_{T^0}^\ell &\coloneqq \{ v \in V_0 \with \operatorname{supp}(v) \subset \mathsf{N}_\ell(T^0) \}, \label{eq:VT0ell}\\
M_{T^0}^\ell &\coloneqq \left\{ \mu \in \mathbb{R}^N \with \mu_{T^c} = 0 \text{ for all } T^c \subset \Omega^c \setminus \mathsf{N}_\ell^c(T^0),\; c \in \{0,1\} \right\}.
\end{align}
\\
\begin{figure}[h!] \centering
\begin{subfigure}{0.49\textwidth} \centering
	\includegraphics[scale=.2]{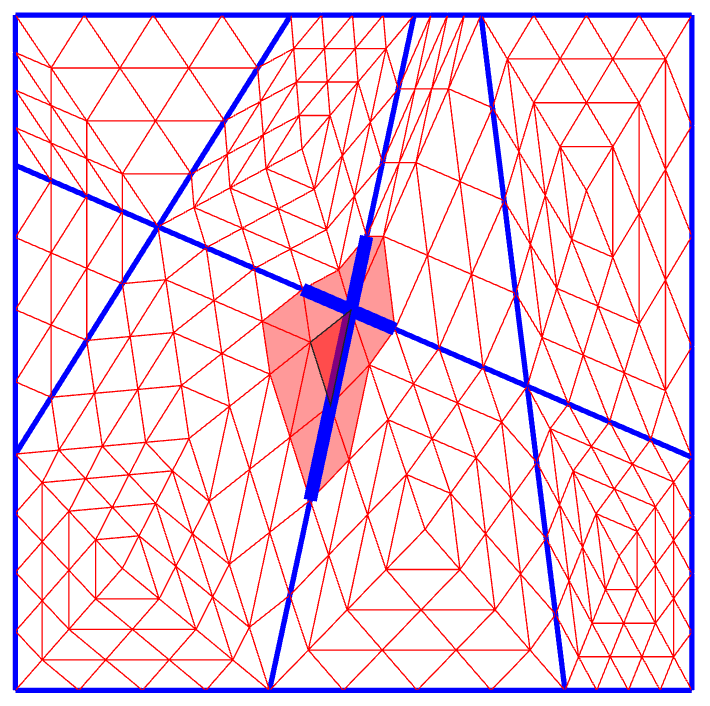}
\end{subfigure}
\caption{Example of a first-order  patch $\mathsf{N}_1(T^0)$, around a coarse bulk element. The shaded red area marks the bulk elements of the patch whereas the opaque red marks the element $T$. The thickened blue lines mark the interface elements belonging to the patch.}
\end{figure}
Localized element correctors can then be defined for any \( v \in V_0 \) as the unique solution \( (\mathcal{K}_{T^0}^\ell v, \lambda^\ell) \in V_{T^0}^\ell \times M_{T^0}^\ell \) satisfying
\begin{subequations}
\label{eq:defKTloc} 
\begin{align}
	& a (\mathcal K_{T^0}^\ell v, w)& +&  &b(w,\lambda^\ell) & &=\quad  &a_{T^0}(\mathcal I_H v,w), &&\text{ for all }w \in V_{T^0}^\ell, && \label{eq:defKTloc1}\\
	& b(\mathcal K_{T^0}^\ell v,\mu)                   &   &         &    & &=\quad  &-b_{T^0}(v-\mathcal I_H v,\mu), &&\text{ for all }\mu \in  M_{T^0}^\ell.&& \label{eq:defKTloc2}
\end{align}
\end{subequations}
We define a localized version of the operator \(\mathcal{R}\) as:
\begin{equation}
\label{eq:defRl}
\mathcal R^\ell v\coloneqq (\mathcal I_H - \mathcal K^\ell)v,\qquad \mathcal K^\ell v \coloneqq \sum_{T^0 \in \mathcal{T}_H^0} \mathcal K_{T^0}^\ell v.
\end{equation}

In what follows, we show that \(\mathcal{R}^\ell\) provides an exponentially accurate \mbox{(quasi-)local} approximation of \(\mathcal{R}\). The proof relies on a localized version of the decay result stated in \cref{thm:dec}, which is formulated in the following corollary.

\begin{corollary}[Local version of exponential decay]\label{thm:decloc}
For all $T^0 \in \TH$  and all $m,\ell \in \mathbb N$, it holds that
\begin{equation*}
	\|\mathcal K^\ell_{T^0} v\|_{V(\Omega\setminus \Nb_{m}(T^0))} \lesssim \exp(-c m)\|\mathcal K^\ell_{T^0} v\|_{V( \mathsf N_\ell (T^0))},\quad v \in V^0,
\end{equation*}
where $c$ is the constant from \cref{thm:dec}.
\end{corollary}
\begin{proof}
We apply \cref{thm:dec} with \( \mathcal{K}_{T^0} \), \( V \), and \( \mathsf{N}_\ell(T^0) \) replaced by \( \mathcal{K}_{T^0}^\ell \), \( V_{T^0}^\ell \), and~\( \mathsf{N}_m(T^0) \), respectively. The proof carries over without modification.
\end{proof}

We now proceed to the exponential approximation result.
\begin{theorem}[Localization error]\label{thm:localizationerror}
For all $\ell \in \mathbb N$, it holds that 
\begin{equation}
	\label{eq:locerrR}
	\|(\mathcal R - \mathcal R^\ell) v\|_{V} \lesssim \ell^{(d-1)/2}\exp(-c \ell)\| v\|_{V},\quad v \in V_0,
\end{equation}
where $c$ is the constant from \cref{thm:dec}.
\end{theorem}
\begin{proof}
Let \( e \coloneqq (\mathcal{R} - \mathcal{R}^\ell)v \) denote the localization error. Using the coercivity of the bilinear form \(a\) from \cref{lem:coercivity}, the fact that $a(\tilde V_H,w) = 0$, and the definition of \(\mathcal{R}^\ell\) from \cref{eq:defRl}, we obtain the following representation of the localization error:
\begin{align}\label{eq:spliterror}
	\|e\|_{V}^2 \lesssim  a(e,e)= -a(\mathcal R^\ell v,e) = -a(\mathcal I_Hv,e) + \sum_{T^0 \in \mathcal{T}_H^0} a(\mathcal K_{T^0}^\ell v,e).
\end{align}
We now analyze each term in the sum on the right-hand side individually.
Recalling the definition of the cut-off function from \cref{eq:eta} and testing \cref{eq:defKTloc1} with the test function \( (1 - \eta)e \in V_{T^0}^\ell \), we obtain that
\begin{align*}
	a(\mathcal K_{T^0}^\ell v,e) &= a(\mathcal K_{T^0}^\ell v,(1-\eta)e + \eta e)\\
	&= -b((1-\eta) e,\lambda^\ell) + a_{T^0}(\mathcal I_H v,(1-\eta)e)+a(\mathcal K_{T^0}^\ell v,\eta e)\\&\eqqcolon \Xi_1+\Xi_2+\Xi_3.
\end{align*}
To estimate the term \( \Xi_1 \), we employ arguments analogous to those used in \cref{eq:cstart}, replacing \( \eta \) with \( 1 - \eta \) and \( \varphi \) with \( e \). Note that, not only \( \varphi \), but also \( e \) satisfies the vanishing mean condition \cref{eq:avgvarphi}, which leads to the following estimate:
\begin{align*}
	|\Xi_1| &= \bigg|\sum_{c \in \{0,1\}}\sum_{K^c \subset R_{T^0}^c} \lambda_{K^c}^\ell q_{K^c}((1-\eta)e)\bigg| \\
	&\lesssim \sum_{c \in \{0,1\}}\Big( H^{-d+c+2}\sum_{K^c\subset R_{T^0}^c}  |\lambda_{K^c}^\ell|^2\Big)^{1/2}  \|\nabla^ce^c\|_{L^2(R_{T^0}^c)},
\end{align*}
where $R_{T^0}^c$ is defined as in \cref{eq:ring}.
To estimate the entries of \( \lambda^\ell \in M_{T^0}^\ell \), we follow the same steps as in the proof of \cref{thm:dec}, specifically \cref{eq:sumlamK1,eq:sumlamK0}, yielding
\begin{equation}\label{eq:estlambda}
	\sum_{c \in \{0,1\}}H^{-d+c+2}\sum_{K^c \subset R_{T^0}^c}  |\lambda^\ell_{T^c}|^2
	\lesssim \|\mathcal{K}_{T^0}^\ell v\|^2_{V(\mathsf N(R_{T^0}))}.
\end{equation}
Inserting this into the above estimate for $|\Xi_1|$ gives
$$
|\Xi_1|\lesssim \|\mathcal{K}_{T^0}^\ell v\|_{V(\mathsf N(R_{T^0}))} \|e\|_{V(R_{T^0})}. 
$$

Next, we consider the term \( \Xi_2 \). Noting that \( a_{T^0}(\mathcal{I}_H v, (1 - \eta) e) = a_{T^0}(\mathcal{I}_H v, e) \), it follows that, upon summing over all bulk elements \( T^0 \in \mathcal{T}_H^0 \), the term \( \Xi_2 \) cancels exactly with the term \( -a(\mathcal{I}_H v, e) \) on the right-hand side of \cref{eq:spliterror}.

To estimate \( \Xi_3 \), we use that the intersection of the supports of the functions~\( \mathcal{K}_{T^0}^\ell v \) and \( \eta e \) is contained in \( R_{T^0} \), which yields that
\begin{align*}
	\Xi_3 \lesssim |a(\mathcal{K}_{T^0}^\ell v, \eta e)| 
	\lesssim \|\mathcal{K}_{T^0}^\ell v\|_{V( R_{T^0})} \, \|\eta e\|_{V( R_{T^0})} 
	\lesssim \|\mathcal{K}_{T^0}^\ell v\|_{V( R_{T^0})} \, \|e\|_{V( \mathsf{N}(R_{T^0}))},
\end{align*}
where in the last step we employed  definition~\cref{eq:restanorm}, the trace inequality~\cref{eq:trace}, and the fact that \( e \) satisfies~\cref{eq:avgvarphi}.

Combining the estimates for the terms \( \Xi_1 \), \( \Xi_2 \), and \( \Xi_3 \), and applying \cref{thm:decloc}, along with the finite overlap of the sets \( \{ \mathsf{N}(R_{T^0}) \with T^0 \in \mathcal{T}_H^0 \} \), we obtain
\begin{align}
	\|e\|_{V}^2&\lesssim \sum_{T^0\in \TH^0} \|\mathcal K_{T^0}^\ell v\|_{V(\mathsf N (R_{T^0}))}\|e\|_{V(\mathsf N(R_{T^0}))}\notag\\
	&\lesssim \exp(-c\ell) \sum_{T^0 \in \TH^0}\| \mathcal K_{T^0}^\ell v\|_{V}\|e\|_{V(\mathsf N(R_{T^0}))}\notag\\
	&\lesssim \exp(-c\ell) \sqrt{\sum_{T^0 \in \TH^0} \|\mathcal K_{T^0}^\ell v\|^2_{V}}\sqrt{\sum_{T^0 \in \TH^0} \| e\|^2_{V(\mathsf N(R_{T^0}))}}\notag\\
	&\lesssim \ell^{(d-1)/2}\exp(-c\ell) \sqrt{\sum_{T^0 \in \TH^0} \|\mathcal K_{T^0}^\ell v\|^2_{V}}\| e\|_{V}.\label{eq:almostfinalest} 
\end{align}
It remains to estimate the term \( \| \mathcal{K}_{T^0}^\ell v \|_{V} \). To this end, we test \cref{eq:defKTloc1} with the function \( \mathcal{K}_{T^0}^\ell v \) and apply the identity obtained by testing \cref{eq:defKTloc2} with \( \lambda^\ell \). Using \cref{lem:coercivity}, this yields the following estimate:
\begin{align}
	\label{eq:estKtlv}
	\|\mathcal K_{T^0}^\ell v\|^2_{V}&\lesssim a(\mathcal K_{T^0}^\ell v,\mathcal K_{T^0}^\ell v)= a_{T^0}(\mathcal{I}_H v,\mathcal K_{T^0}^\ell v) + b_{T^0}(v-\mathcal{I}_H v,\lambda^\ell).
\end{align}
We bound the first term on the right-hand side of the latter estimate as
\begin{align*}
	|a_{T^0}(\mathcal{I}_H v,\mathcal K_{T^0}^\ell v)| \lesssim \|\mathcal I_H v\|_{V(T^0)}\|\mathcal K_{T^0}^\ell v\|_{V(T^0)} \lesssim \|v\|_{V(\mathsf N(T^0))}\|\mathcal K_{T^0}^\ell v\|_{V(T^0)},
\end{align*}
where we have used \cref{eq:IH} and the trace inequality, cf.~\cref{eq:trace}. Above $\|\cdot\|_{V(T^0)}$ is used as an abbreviation for $\|\cdot\|_{V(T^0,\partial T^0)}$. 

For the second term, we obtain the estimate
\begin{align*}
	&|b_{T^0}(v-\mathcal I_H v,\lambda^\ell)|\\
	&\quad  \lesssim H^{-d/2+1}|\lambda^\ell_{T^0}|\|\nabla^0 v^0\|_{L^2(T^0)} + \Big(H^{-d+3}\sum_{T^1\subset \partial T^0} |\lambda_{T^1}^\ell|^2\Big)^{1/2}\|\nabla^1 v^1\|_{L^2(\partial T^0)}\\
	&\quad \lesssim \|\mathcal K_{T^0}^\ell v \|_{V(\mathsf N(T^0))}\|v \|_{V(\mathsf N(T^0))},
\end{align*}
where we used arguments similar to those in the proof of \cref{eq:estlambda} to estimate the components of \( \lambda^\ell \). Inserting the two preceding estimates into \cref{eq:estKtlv} yields
\begin{equation*}
	\|\mathcal K_{T^0}^\ell v\|^2_{V(\mathsf N_\ell(T^0))}\lesssim \|v\|_{V(\mathsf N_\ell(T^0))}\|\mathcal K_{T^0}^\ell v\|_{V(\mathsf N_\ell(T^0))},
\end{equation*}
and, in turn, inserting this estimate into \cref{eq:almostfinalest} yields the assertion.
\end{proof}

\section{Practical multiscale method}\label{sec:pracmethod}
\begin{figure}[h!] \centering
\begin{subfigure}{.49\textwidth} \centering
	\includegraphics[scale=.28]{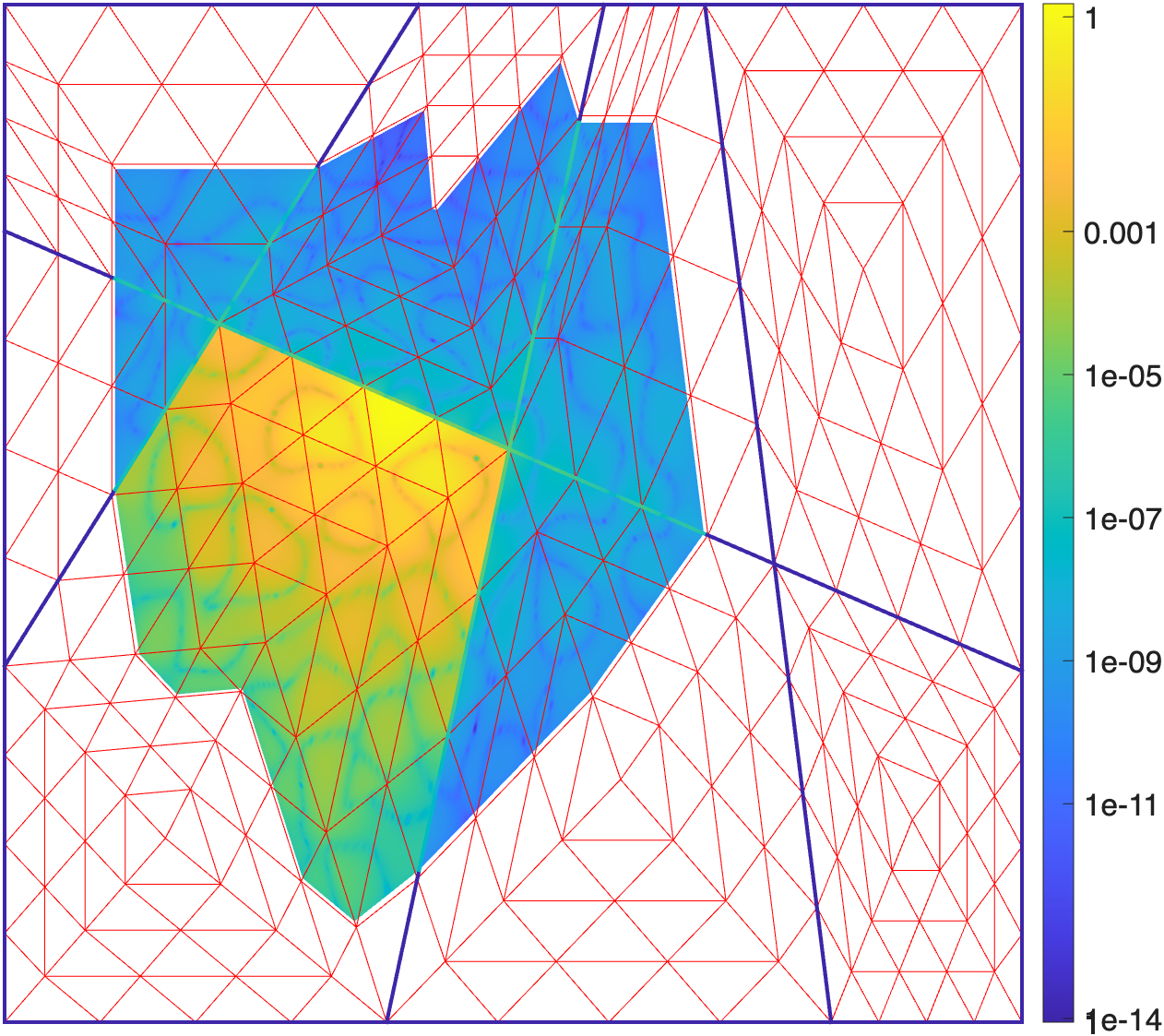}\caption{}\label{subfig:basis_functions_bulk}
\end{subfigure}%
\begin{subfigure}{.49\textwidth} \centering
	\includegraphics[scale=.28]{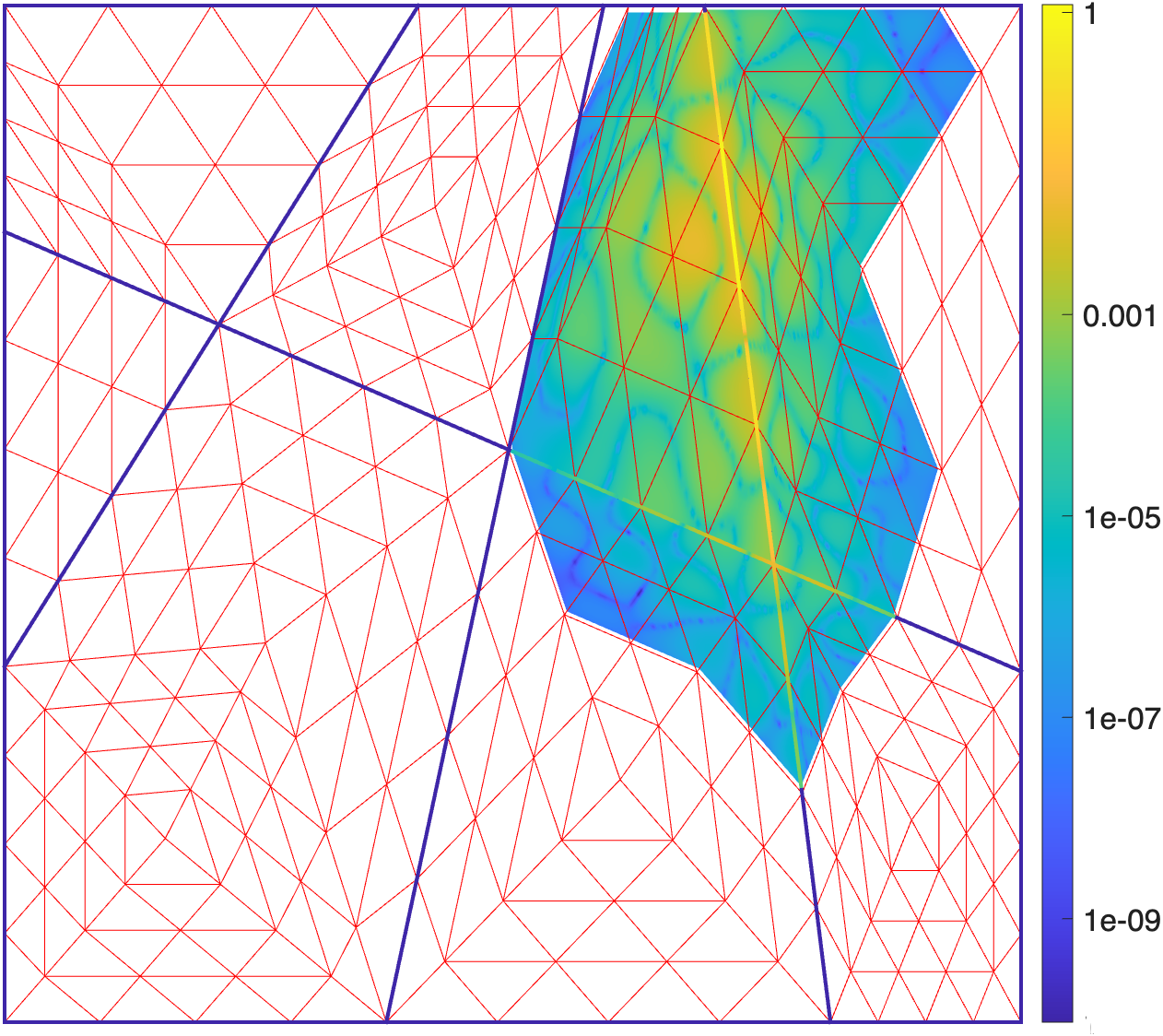}\caption{}\label{subfig:basis_functions_interface}
\end{subfigure}%
\caption{Examples of basis functions corresponding to a coarse mesh bulk element \subref{subfig:basis_functions_bulk} and a coarse mesh interface element \subref{subfig:basis_functions_interface} using a patch spreading $\ell = 4$ layers.}\label{fig:basis_functions}
\end{figure}

In this section, we introduce a practical multiscale method based on locally computable basis functions. This is justified by the exponential decay of the globally defined prototypical basis functions. We define the localized basis functions as
\begin{equation}
\label{eq:locbasis}
\tilde \varphi_{T^c}^\ell \coloneqq \mathcal{R}^\ell b_{T^c},
\end{equation}
and define the localized multiscale space as their span:
\begin{equation*}
\tilde V_H^\ell \coloneqq \operatorname{span}\{\tilde \varphi_{T^c}^\ell\with T^c \in \mathcal{T}_H^c,\, c \in \{0,1\}\}.
\end{equation*}
\\
We emphasize that the operator \( \mathcal{R}^\ell \) depends on its argument only through its QOI, as introduced in \cref{defQOI}. Consequently, the definition \cref{eq:locbasis} of \( \tilde\varphi^\ell_{T^c} \) is independent of the specific choice of bubble functions, provided that the Kronecker-delta condition~\cref{eq:kdbubble} on the QOI is satisfied. 
Note that, for a bulk element \( T^0 \in \mathcal T_H^0 \), computing the associated basis function \( \tilde \varphi_{T^0}^\ell \) requires solving local problems of the form \cref{eq:defKTloc} for all bulk elements in the first-order patch surrounding \( T^0 \). Similarly, for an interface element \( T^1 \in \mathcal T_H^1 \), the basis function \( \tilde \varphi_{T^1}^\ell \) is obtained by solving local problems of the form \cref{eq:defKTloc} for each bulk element that shares the face \( T^1 \). Figure~\ref{fig:basis_functions} shows examples of a bulk basis function \( \tilde \varphi_{T^0}^\ell \) and an interface basis function \( \tilde \varphi_{T^1}^\ell \) for \(\ell = 4\).
\\

The LOD method seeks the unique function $\tilde u_H^\ell\in \tilde V_H^\ell$ such that
\begin{align}\label{eq:locmethod}
a(\tilde u_H^\ell,\tilde v_H^\ell)=F(\tilde v_H^\ell), 
\quad \text{for all } \tilde v_H^\ell \in \tilde V_H^\ell.
\end{align}
The following theorem provides a convergence result for the LOD approximation.

\begin{theorem}[Localized method]
\label{thm:convergenceloc}
The localized multiscale method~\cref{eq:locmethod} is well-posed, and, for any right-hand side \( F \) of the form \cref{eq:rhs}, with components \( f^0 \in H^s(\Omega^0) \) and \( f^1 \in H^s(\Omega^1) \), where \( s \in \{0,1\} \), the following error estimate holds:
\begin{align}
	\|u - \tilde{u}_H^{\ell}\|_V & \lesssim H^{1+s} |f|_s + \ell^{(d-1)/2} \exp(-c \ell) |f|_0, \label{eq:errestpracH1}
\end{align}
where the (semi-)norm \( |\cdot|_s \) is defined in \cref{eq:fnorm}.
\end{theorem}
\begin{proof}
After establishing Theorem \ref{thm:convergenceprot} and Theorem \ref{thm:localizationerror} the proof follows the arguments of \cite[Thm.~5.2]{MalP20} closely. We omit the details for the sake of brevity.
\end{proof}

We emphasize that, for a practical computation of the LOD basis functions defined in \cref{eq:locbasis}, the problems \cref{eq:defKTloc}, which are local but infinite-dimensional, still need to be solved. In practice, we approximate the solutions to these problems using a (local) fine-scale discretization. For codimensions $c \in \{0,1\}$, let $\mathcal{T}_h^c$ denote a mesh obtained by uniform refinement of the coarse mesh $\mathcal{T}_H^c$. The fine mesh size $h<H$ should be chosen such that $\mathcal{T}_h^c$ resolves the fine-scale oscillations of the coefficients of the problem at hand.
Throughout the derivation of the method, we replace the space $V_0$ by its finite-dimensional subspace $V_{0,h}$, consisting of continuous, piecewise linear functions defined with respect to $\mathcal{T}_h^c$. Specifically, for problem \cref{eq:defKTloc}, this means posing it on a finite-dimensional version of $V_{T^0}^\ell$ as defined in \cref{eq:VT0ell}, but with~$V_0$ replaced by $V_{0,h}$.
In the fully discrete convergence analysis, $V_0$ is likewise replaced by $V_{0,h}$, and most arguments carry over directly; see, e.g., \cite[Ch.~4.4]{MalP20}, for a discussion of technical details that also arise in the present setting. This yields an a priori error estimate for the fully discrete LOD approximation similar to \cref{thm:convergenceloc}, but with respect to the fine-scale finite element solution in $V_{0,h}$. An estimate with respect to the weak solution of the original PDE then follows from the triangle inequality combined with standard finite element approximation results; see, e.g., \cite[Thm.~2]{Hellman2023} for the present mixed-dimensional setting.

\section{Generalization to complex-shaped coarse elements}\label{sec:complexelements}

The assumption made in \cref{sec:proto} and used throughout, namely that the coarse meshes $\mathcal{T}_H^c$, $c \in \{0,1\}$, are simplicial or quadrilateral/hexahedral, can be relaxed. A closer inspection reveals that the construction of the method and its error analysis relies only on the ability to:
\begin{enumerate}
\item construct bubble functions satisfying the Kronecker--delta property~\eqref{eq:kdbubble} and the inverse inequality~\eqref{eq:stabbubble};\label{prop1}
\item construct a suitable partition of unity on the meshes $\mathcal{T}_H^c$, $c \in \{0,1\}$, for defining a quasi-interpolation operator in place of definition~\eqref{eq:IH};\label{prop2}

\item establish local Poincar\'e-type inequalities used to prove that properties~\eqref{eq:propIH} hold for the resulting partition-of-unity-based quasi-interpolation.\label{prop3}
\end{enumerate}

\begin{figure}[h!] \centering
\begin{subfigure}{.475 \textwidth} \centering
	\includegraphics[scale=.2]{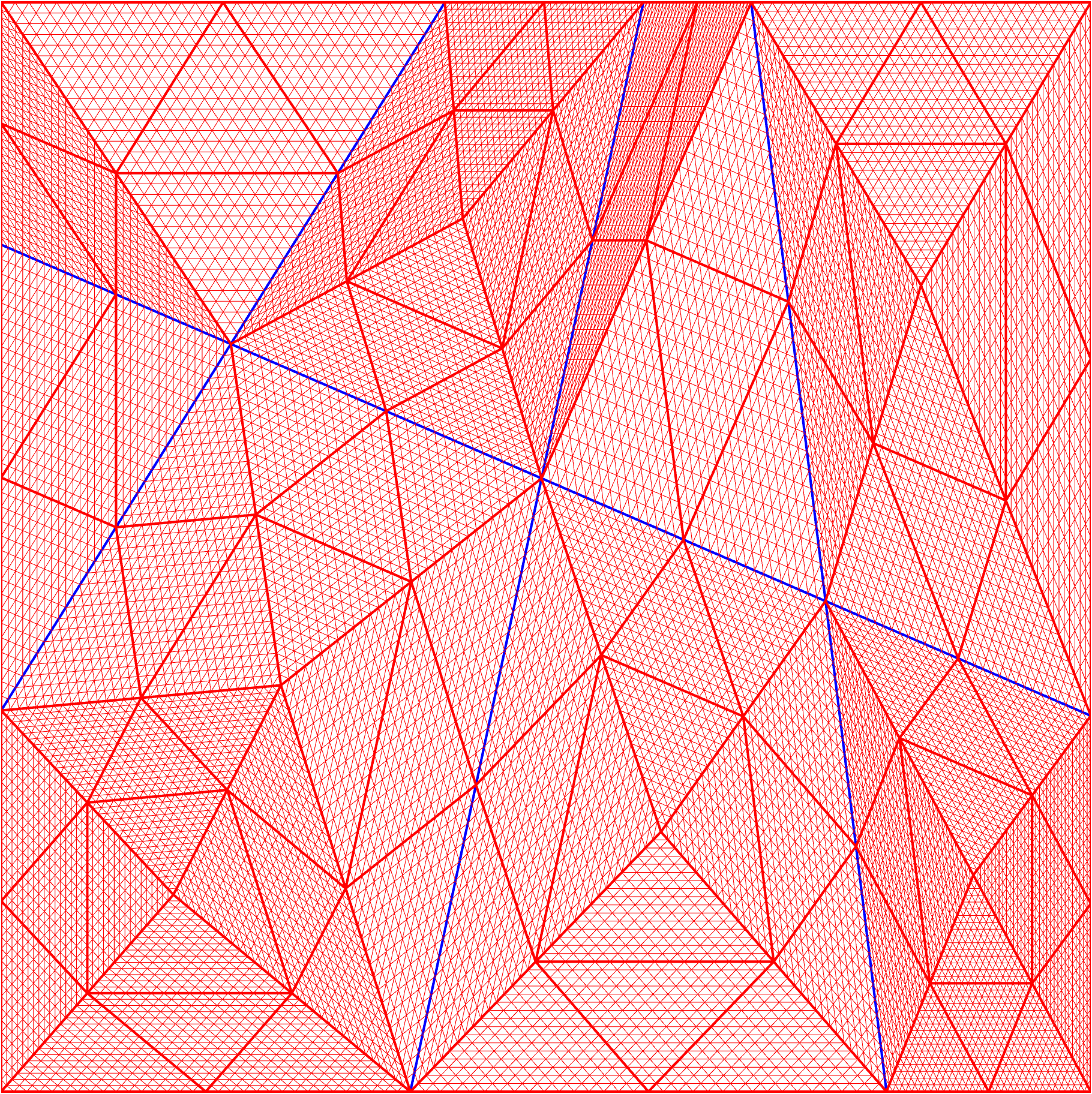}
	\caption{}\label{fig:complex_shaped_domain_orig}
\end{subfigure}
\hfill
\begin{subfigure}{.475 \textwidth} \centering
	\includegraphics[scale=.2]{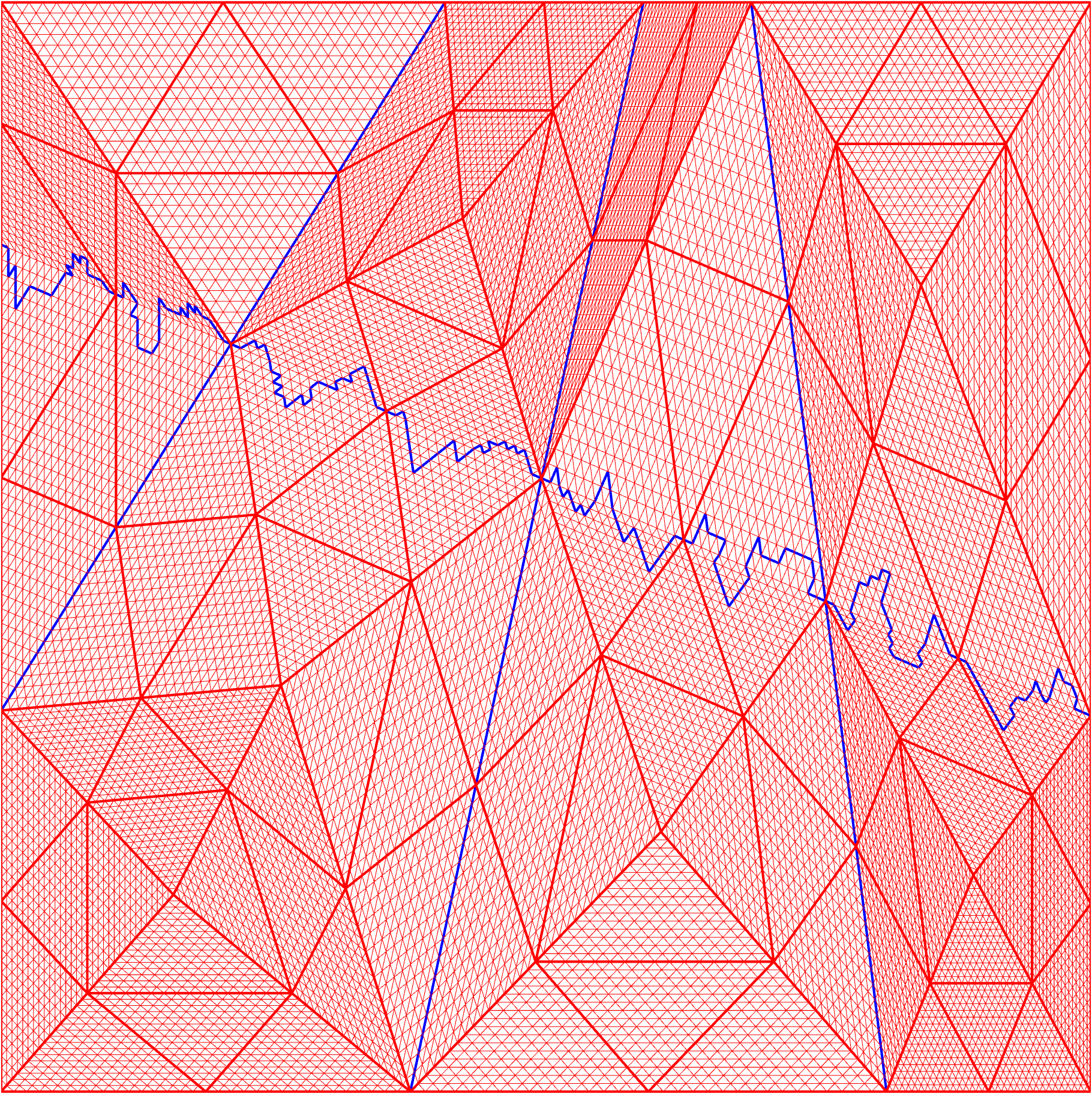}
	\caption{}\label{fig:complex_shaped_domain_new}
\end{subfigure}
\caption{A domain with straight interfaces (blue lines) and simplicial coarse mesh elements (red thick lines) is shown in \subref{fig:complex_shaped_domain_orig}. A modified domain with an interface that follows the fine mesh (red thin lines) and leads to non-convex, complex-shaped coarse mesh elements is shown in \subref{fig:complex_shaped_domain_new}.}\label{fig:complex_shaped_domain}
\end{figure}

Such assumptions also hold in a more general setting for shape-regular and quasi-uniform coarse meshes $\mathcal{T}_H^c$, $c \in \{0,1\}$, consisting of simply connected elements formed as unions of fine mesh elements from $\mathcal{T}_h^c$. See Figure~\ref{fig:complex_shaped_domain} for an example of such a mesh. For these meshes, quasi-uniformity and non-degeneracy can be characterized as follows: for each element $T^c \in \mathcal{T}_H^c$, with $c \in \{0,1\}$, there exist inscribed and circumscribed balls $B_R^c(x_0) \subset T^c \subset B_{R'}^c(x_0)$ such that
\[
\rho_0 H \leq R \leq R' \leq \rho_1 H,
\]
where $\rho_0, \rho_1 > 0$ are moderate constants.

\cref{prop1} holds in the above general setting since a bubble function can be constructed as \( v_{T^c}(x) = \alpha \max(0, R - |x - x_0|) \) with support on \( T^c \). By choosing the constant \( \alpha \) appropriately, this function satisfies the Kronecker--delta property~\eqref{eq:kdbubble}. Clearly, also the inverse inequality~\eqref{eq:stabbubble} is satisfied.

To construct a partition of unity suitable for defining a quasi-interpolation operator replacing \cref{eq:IH}, as required in \cref{prop2}, we first introduce open covers of \(\Omega^c\) for \(c \in \{0,1\}\). The subdomains in these covers correspond to elements \(T^c \in \mathcal{T}_H^c\) and are defined as the first-order element patches with respect to \(\mathcal{T}_H^c\), i.e., $U_{T^c} = \mathsf{N}^c(T^c)$ with patches defined as in~\cref{eq:patchesc}.
Define the functions
\[
\tilde{\Lambda}_{T^c}(x) = \begin{cases}
\mathrm{dist}^c(x, \partial U_{T^c}), & x \in U_{T^c}, \\
0, & \text{otherwise},
\end{cases}
\]
where \(\mathrm{dist}^c(x, \partial U_{T^c})\) denotes the intrinsic distance from \(x\) to the boundary \(\partial U_{T^c}\), considering only paths contained entirely in $\Omega^c$. Using these functions, the partition of unity for \( c \in \{0,1\} \) consists of functions \(\Lambda_{T^c}\) defined for each \( T^c \in \mathcal{T}_H^c \) by
\begin{equation}\label{eq:PofU}
\Lambda_{T^c}(x) = \frac{\tilde{\Lambda}_{T^c}(x)}{\sum_{S^c \in \mathcal{T}_H^c} \tilde{\Lambda}_{S^c}(x)}.
\end{equation}
Since the coarse meshes are assumed to be quasi-uniform and shape-regular, the supports of the partition of unity functions overlap sufficiently to guarantee that the denominator in \eqref{eq:PofU} remains uniformly positive. Moreover, it follows that each~\(\Lambda_{T^c}\) is Lipschitz continuous, with a Lipschitz constant proportional to \(H^{-1}\). A quasi-interpolation operator can then be constructed as
\begin{equation*}
\mathcal I_H^c v^c\coloneqq \sum_{T^c \in \mathcal T_H^c} q_{T^c}(v^c) \Lambda_{T^c};
\end{equation*}
see \cite{Car99} for the analysis of such operators, and also \cite{HMM23}.

Finally, concerning \cref{prop3}, we assume that each element \( T^c \in \mathcal{T}_H^c \) satisfies a local Poincaré inequality on its first-order patch:
\begin{equation}\label{eq:Poincarepatch}
\|v - q_{\mathsf{N}^c(T^c)}(v)\|_{L^2(\mathsf{N}^c(T^c))} \leq C H \|\nabla^c v\|_{L^2(\mathsf{N}^c(T^c))},
\end{equation}
for all \( v \in H^1(\mathsf{N}^c(T^c)) \), where \( q_{\mathsf{N}^c(T^c)} \) is defined as in \cref{defQOI} but with \( \mathsf{N}^c(T^c) \) in place of \( T^c \).
According to \cite{Ruiz12}, for an inequality of the form \eqref{eq:Poincarepatch} to hold, three conditions are required: a uniform bound on the size of \( \mathsf{N}^c(T^c) \); the existence, for each \( \mathsf{N}^c(T^c) \), of a fixed finite cone \( C \) such that every boundary point of \( \mathsf{N}^c(T^c) \) is the vertex of a cone \( C_x \), congruent to \( C \) and contained in \( \mathsf{N}^c(T^c) \); and a constant \(\delta_0 > 0\) such that for all \(\delta \in (0, \delta_0)\), the inner set
\[
\mathsf{N}^c_\delta(T^c) \coloneqq \{ x \in \mathsf{N}^c(T^c) \with \mathrm{dist}^c(x, \partial \mathsf{N}^c(T^c)) > \delta \}
\]
is connected. These assumptions exclude domains with narrow regions but are not very restrictive in our context, as our elements are constructed as unions of connected fine-scale elements.

\section{Numerical experiments}\label{sec:numerics}
In this section, we present a numerical investigation of the proposed multiscale method. The first two numerical examples examine the convergence of the energy norm of the error, $\|u_h-\tilde{u}_{H,h}^\ell \|_a$, with respect to the coarse mesh size $H$ for the localized method defined in \cref{eq:locmethod}. Here, $u_h$ denotes the solution obtained from the fitted finite element method on the mesh $\mathcal{T}_h$, as introduced in \cite{Hellman2023}, and $\tilde{u}_{H,h}^\ell$ denotes the solution to the discretized version of problem \cref{eq:locmethod}. The first example considers a problem with smoothly varying coefficients, while the second addresses rapidly oscillating coefficients $A^0$ and $A^1$. The third example investigates the localization error of the method, measured in the energy norm, with respect to the oversampling parameter $\ell$, for a coarse mesh containing complex-shaped elements. In this case, a simplified variant of  method~\cref{eq:locmethod} is employed. Finally, in the last example we use a more complex domain and an example that highlights the multiscale features more clearly. In this example, we look at the convergence following $\ell = \log_2(1/H)$ in order to numerically verify the theoretical results from \cref{thm:convergenceloc}. Throughout this section, $H$ and $h$ denote the coarse and fine mesh sizes, respectively.

\begin{figure}[h!] \centering
\includegraphics{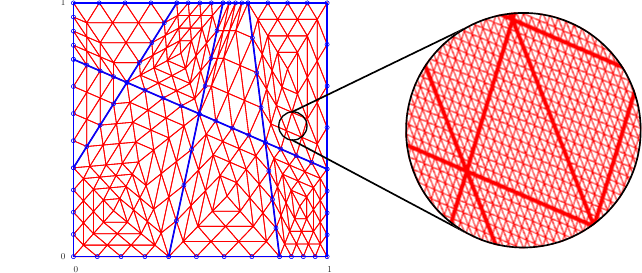}
\caption{The domain $\Omega$ used for the first two numerical experiments. Here, the blue lines are the interfaces $\Omega^1$ and the coarse mesh corresponding to $H = \tfrac14$ is shown in red. The fine-scale mesh is shown in the zoomed-in part.}\label{fig:num_domain}
\end{figure}

We consider the domain $\Omega = (0,1)^2$ with four interfaces, as illustrated in \cref{fig:num_domain}. The convergence with respect to the coarse mesh size $H$ is studied for $H \in \{1, \tfrac{1}{2}, \tfrac{1}{4}, \tfrac{1}{8}, \tfrac{1}{16}\}$. The corresponding coarse meshes are generated by successive uniform refinement of an initial coarse mesh $\mathcal{T}_H$ with $H = 1$, capturing the geometry. The mesh $\mathcal{T}_h$, with mesh size $h = 1/64$, is obtained by uniform refinement of the finest coarse mesh corresponding to $H = \tfrac{1}{16}$. \Cref{fig:num_domain} shows the coarse mesh for $H = \tfrac{1}{4}$, along with a zoomed-in view of the underlying fine mesh. Note that the mesh sizes $H$ and $h$ have to be understood relative to the subdomain size.

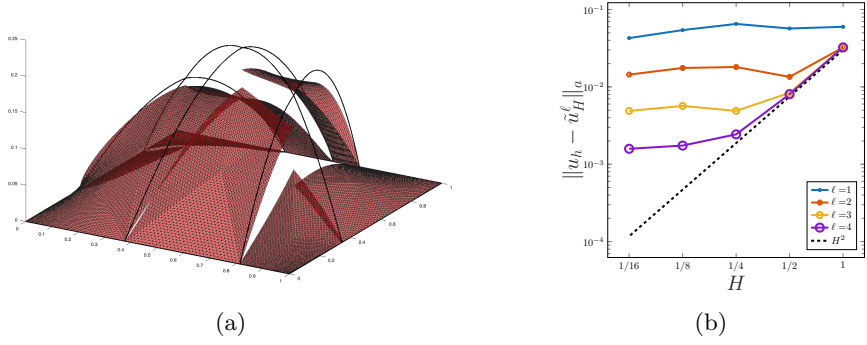
\begin{figure}[h!]
\centering
\begin{subfigure}{.49 \textwidth} \centering
	\includegraphics[width=\textwidth]{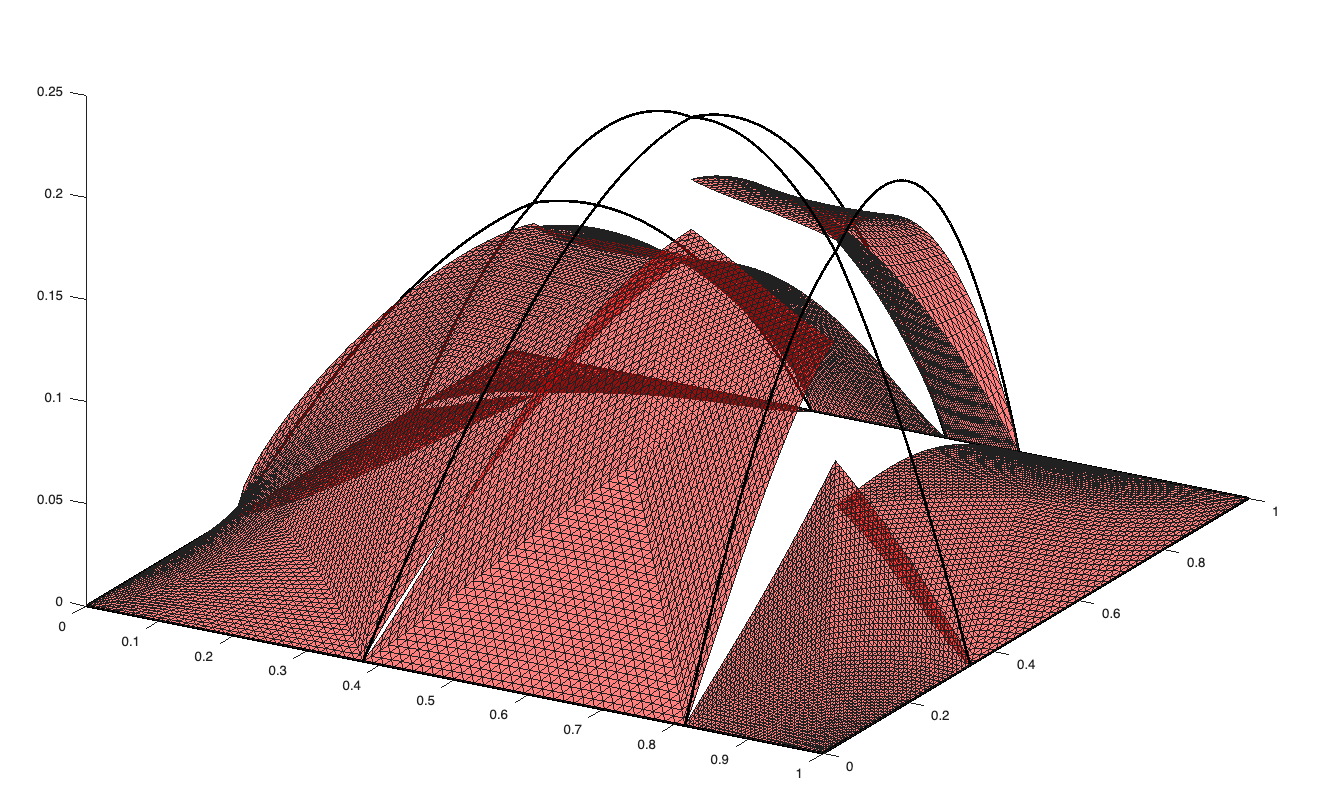}
	\caption{}\label{fig:conv_Ai1_Aj1_fi_sin_fj_2x_sol}
\end{subfigure}
\begin{subfigure}{.49 \textwidth} \centering
        %
        %
        \definecolor{mycolor1}{rgb}{0.06600,0.44300,0.74500}%
        \definecolor{mycolor2}{rgb}{0.86600,0.32900,0.00000}%
        \definecolor{mycolor3}{rgb}{0.92900,0.69400,0.12500}%
        \definecolor{mycolor4}{rgb}{0.52100,0.08600,0.81900}%
        \definecolor{mycolor5}{rgb}{0.12941,0.12941,0.12941}%
        \begin{tikzpicture}[scale=.38]
        
        \begin{axis}[%
        width=3.5in,
        height=3.5in,
        at={(0.948in,0.802in)},
        scale only axis,
        xmode=log,
        xtick={0.0625,0.125,0.25,0.5,1},
        xticklabels={{1/16},{1/8},{1/4},{1/2},{1}},
        xminorticks=true,
        xlabel style={font=\color{mycolor5}\huge, align=center},
        xlabel={$H$},
        ymode=log,
        yminorticks=true,
        ylabel style={font=\color{mycolor5}\huge},
        ylabel={$\|u_h - \tilde u_{H,h}^\ell\|_a$},
        axis background/.style={fill=white},
        legend style={at={(0.97,0.03)}, anchor=south east, legend cell align=left, align=left,font=\huge}
        ]
        \addplot [color=mycolor1, mark=o, line width = 2, mark size = 1, mark options={solid, mycolor1}]
          table[row sep=crcr]{%
        1	0.059842420607469\\
        0.5	0.056901293431663\\
        0.25	0.065162833314657\\
        0.125	0.054244565367396\\
        0.0625	0.042774906725813\\
        };
        \addlegendentry{$\ell = $1}
        
        \addplot [color=mycolor2, mark=o, line width = 2,  mark size = 2, mark options={solid, mycolor2}]
          table[row sep=crcr]{%
        1	0.032383192571294\\
        0.5	0.013470390061628\\
        0.25	0.018070212548787\\
        0.125	0.017575441946848\\
        0.0625	0.014424467297749\\
        };
        \addlegendentry{$\ell = $2}
        
        \addplot [color=mycolor3, mark=o, line width = 2,  mark size = 3, mark options={solid, mycolor3}]
          table[row sep=crcr]{%
        1	0.032433647860598\\
        0.5	0.008380031515889\\
        0.25	0.004896393971657\\
        0.125	0.005668379081384\\
        0.0625	0.004885013943305\\
        };
        \addlegendentry{$\ell = $3}
        
        \addplot [color=mycolor4, mark=o, line width = 2,  mark size = 4, mark options={solid, mycolor4}]
          table[row sep=crcr]{%
        1	0.032433481153816\\
        0.5	0.008044422399445\\
        0.25	0.002434228247662\\
        0.125	0.001743181365697\\
        0.0625	0.001583813999097\\
        };
        \addlegendentry{$\ell = $4}
        
        \addplot [color=black, line width = 2, dashed]
          table[row sep=crcr]{%
        1	0.03\\
        0.5	0.0075\\
        0.25	0.001875\\
        0.125	0.00046875\\
        0.0625	0.0001171875\\
        };
        \addlegendentry{     $H^2$}
        
        \end{axis}
        \end{tikzpicture}
	\caption{}\label{fig:conv_Ai1_Aj1_fi_sin_fj_2x_conv}
\end{subfigure}
\caption{In \subref{fig:conv_Ai1_Aj1_fi_sin_fj_2x_sol}, the solution to the problem defined by $A^0 = 1$, $A^1=1$, $B^1=1$, $f^0 = \sin(\pi x) \sin(\pi y)$, and $f^1= x + 2y$ is shown. In \subref{fig:conv_Ai1_Aj1_fi_sin_fj_2x_conv}, we see the convergence results, the energy norm of the error between the reference solution $u_h$ and approximations $\tilde u_{H,h}^\ell$ with respect to $H$ is shown for different oversampling parameters $\ell$.}\label{fig:conv_Ai1_Aj1_fi_sin_fj_2x}
\end{figure}

In the first numerical experiment, we solve the problem \eqref{eq:weak_problem} on $\Omega$ with data
$A^0 = 1$, $A^1 = 1$, $B^1 = 1$, $f^0 = \sin(\pi x)\sin(\pi y)$, and $f^1 = x + 2y$, and plot the resulting errors for fixed oversampling parameters $\ell$, as a function of the coarse mesh size $H$. 
\cref{fig:conv_Ai1_Aj1_fi_sin_fj_2x}\subref{fig:conv_Ai1_Aj1_fi_sin_fj_2x_sol} shows the finite element reference solution, while
\subref{fig:conv_Ai1_Aj1_fi_sin_fj_2x_conv} shows the convergence behavior of the proposed method.
The error is measured in the energy norm relative to the finite element reference solution.
We observe that, for fixed $\ell$, the error first decreases with optimal order $\mathcal O(H^2)$. Beyond a certain point, the localization error becomes dominant, and the convergence curve reaches a plateau. The level of this plateau decreases with increasing oversampling parameter $\ell$. This behavior is consistent with the theoretical result in \cref{thm:convergenceloc}.
Importantly, for fixed $\ell$, the error does not increase again after reaching the plateau, an undesirable effect observed, for example, in \cite{MalP14,Maier2021}. These improved stability properties of the proposed method can be attributed to the refined localization procedure introduced in \cref{sec:localization}.

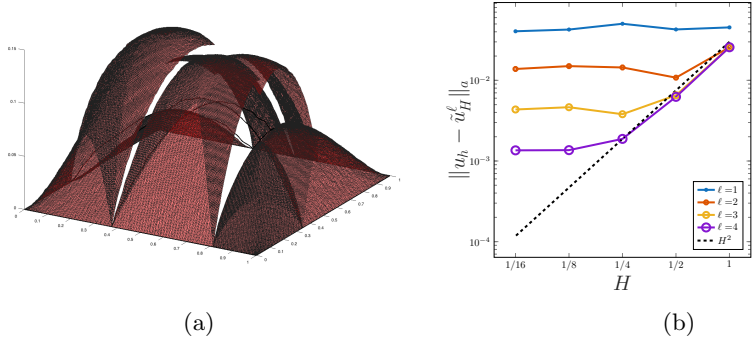
\begin{figure}[h!] \centering
\begin{subfigure}{.49 \textwidth} \centering
	\includegraphics[width=\textwidth]{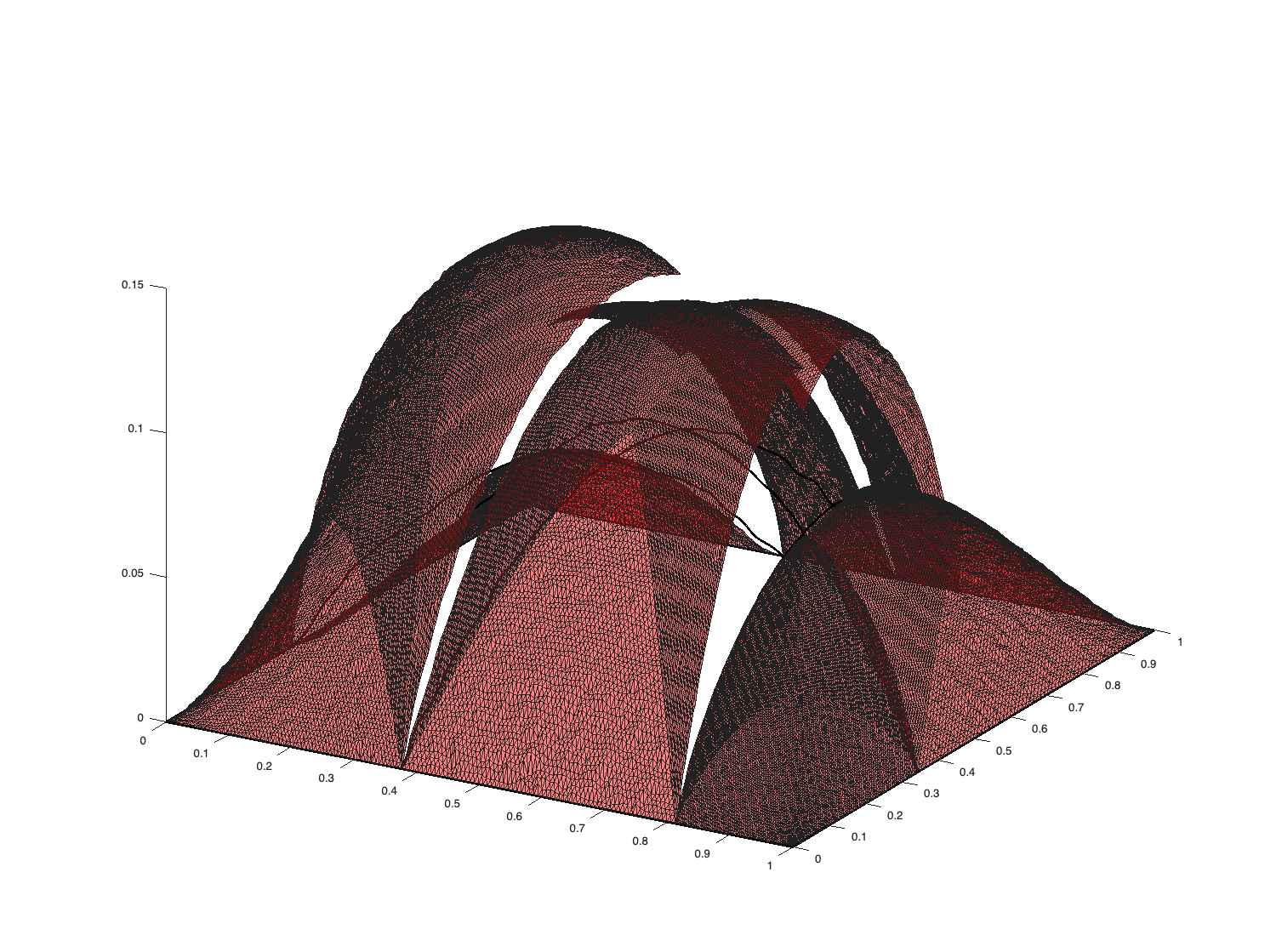}
	\caption{}\label{fig:conv_randomAi_sol}
\end{subfigure}
\begin{subfigure}{.49 \textwidth} \centering
        %
        %
        \definecolor{mycolor1}{rgb}{0.06600,0.44300,0.74500}%
        \definecolor{mycolor2}{rgb}{0.86600,0.32900,0.00000}%
        \definecolor{mycolor3}{rgb}{0.92900,0.69400,0.12500}%
        \definecolor{mycolor4}{rgb}{0.52100,0.08600,0.81900}%
        \definecolor{mycolor5}{rgb}{0.12941,0.12941,0.12941}%
        \begin{tikzpicture}[scale=.38]
        
        \begin{axis}[%
        	width=3.5in,
        	height=3.5in,
        	at={(0.948in,0.802in)},
        	scale only axis,
        	xmode=log,
        	xtick={0.0625,0.125,0.25,0.5,1},
        	xticklabels={{1/16},{1/8},{1/4},{1/2},{1}},
        	xminorticks=true,
        	xlabel style={font=\color{mycolor5}\huge, align=center},
        	xlabel={$H$},
        	ymode=log,
        	yminorticks=true,
        	ylabel style={font=\color{mycolor5}\huge},
        	ylabel={$\|u_h - \tilde u_{H,h}^\ell\|_a$},
        	axis background/.style={fill=white},
        	legend style={at={(0.97,0.03)}, anchor=south east, legend cell align=left, align=left,font=\huge}
        ]
        \addplot [color=mycolor1, mark=o, line width = 2, mark size = 1, mark options={solid, mycolor1}]
          table[row sep=crcr]{%
        1	0.045335547395428\\
        0.5	0.042767106292551\\
        0.25	0.050279626830323\\
        0.125	0.042627873241417\\
        0.0625	0.040583937600623\\
        };
        \addlegendentry{$\ell = $1}
        
        \addplot [color=mycolor2, mark=o, line width = 2, mark size = 2, mark options={solid, mycolor2}]
          table[row sep=crcr]{%
        1	0.025683444007209\\
        0.5	0.01079433560073\\
        0.25	0.014441572198531\\
        0.125	0.01500538361595\\
        0.0625	0.013828634550154\\
        };
        \addlegendentry{$\ell = $2}
        
        \addplot [color=mycolor3, mark=o, line width = 2, mark size = 3, mark options={solid, mycolor3}]
          table[row sep=crcr]{%
        1	0.025712709054594\\
        0.5	0.006531002214135\\
        0.25	0.003802192779816\\
        0.125	0.004642232883428\\
        0.0625	0.004350696801523\\
        };
        \addlegendentry{$\ell = $3}
        
        \addplot [color=mycolor4, mark=o, line width = 2, mark size = 4, mark options={solid, mycolor4}]
          table[row sep=crcr]{%
        1	0.025712443151336\\
        0.5	0.006242434650653\\
        0.25	0.001879541291834\\
        0.125	0.001363050747599\\
        0.0625	0.001355202781147\\
        };
        \addlegendentry{$\ell = $4}
        
        \addplot [color=black, line width = 2, dashed]
        table[row sep=crcr]{%
        	1	0.03\\
        	0.5	0.0075\\
        	0.25	0.001875\\
        	0.125	0.00046875\\
        	0.0625	0.0001171875\\
        };
        \addlegendentry{     $H^2$}
        \end{axis}
        \end{tikzpicture}%
	\caption{}\label{fig:conv_randomAi_conv}
\end{subfigure}
\caption{In \subref{fig:conv_randomAi_sol}, the solution of the problem with random, piecewise constant $A^0$ on the fine mesh, with random values in the range $[0.01,~1]$. Here, $A^1=~\sin(30 \pi x) \sin(30 \pi y) + 2$, $B^1=1$, $f^0 = \sin(\pi x) \sin(\pi y)$, and $f^1= x + 2y$. In \subref{fig:conv_randomAi_conv}, the convergence of said problem with regards to the coarse mesh size $H$, for different oversampling parameters $\ell$ is shown. }\label{fig:conv_randomAi}
\end{figure}

In the second numerical example, we consider coefficients $A^0$ and $A^1$ that vary rapidly on the fine scale. The coefficient $A^0$ is chosen to be piecewise constant on the fine mesh $\mathcal T_h$, taking random values in the range $[0.01, 1]$, while $A^1$ is defined as $A^1 =~\sin(30\pi x)\sin(30\pi y) + 2$. For this example, again we let $B^1 = 1$, $f^0~=~\sin(\pi x) \sin(\pi y)$, and $f^1= x + 2y$. \cref{fig:conv_randomAi}\subref{fig:conv_randomAi_sol} shows the corresponding finite element solution, and \subref{fig:conv_randomAi_conv} presents the resulting convergence behavior. The observed convergence is qualitatively similar to that of the previous numerical example.

\begin{figure}[h!] \centering
\begin{subfigure}{.49 \textwidth} \centering
	\includegraphics[width=\textwidth]{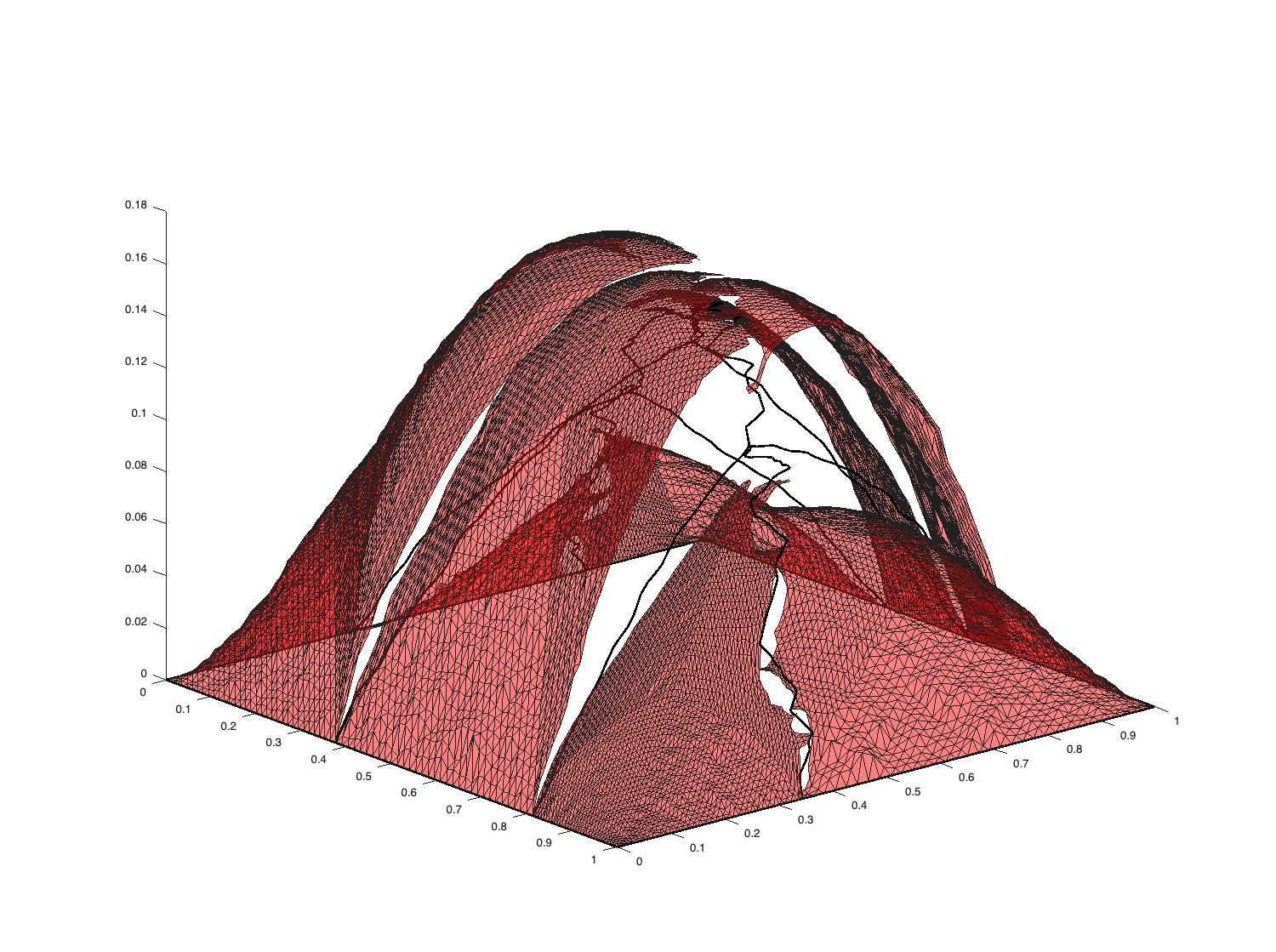}
	\caption{}\label{fig:complex_shaped_sol}
\end{subfigure}
\begin{subfigure}{.49 \textwidth} \centering
        %
        %
        \definecolor{mycolor1}{rgb}{0.06600,0.44300,0.74500}%
        \definecolor{mycolor2}{rgb}{0.12941,0.12941,0.12941}%
        \begin{tikzpicture}[scale=.38]
        
        \begin{axis}[%
        width=3.5in,
        height=3.5in,
        at={(0.948in,0.802in)},
        scale only axis,
        xtick={1, 2, 3, 4, 5, 6},
        xlabel style={font=\color{mycolor2}\huge},
        xlabel={$\ell$},
        ymode=log,
        yminorticks=true,
        ylabel style={font=\color{mycolor2}\huge},
        	ylabel={$\|u_h - \tilde u_{H,h}^\ell\|_a$},
        legend style={font=\huge},
        axis background/.style={fill=white}
        ]
        \addplot [color=mycolor1, mark size = 3, line width = 2, mark=o, mark options={solid, mycolor1}]
          table[row sep=crcr]{%
        1	0.524951357613358\\
        2	0.126788150645751\\
        3	0.019701416608526\\
        4	0.002568974026859\\
        5	0.000307441986633\\
        6	3.7131984985e-05\\
        };
        \addlegendentry{$H = \tfrac12$};
        \addplot [color=black, dotted, line width = 2]
        table[row sep=crcr]{%
          1.000000000000000   0.135335283236613\\
        2.000000000000000   0.018315638888734\\
        3.000000000000000   0.002478752176666\\
        4.000000000000000   0.000335462627903\\
        5.000000000000000   0.000045399929762\\
        6.000000000000000   0.000006144212353\\
        };
        \addlegendentry{$\exp(-2\ell)$};
        \end{axis}
        \end{tikzpicture}%
	\caption{}\label{fig:complex_shaped_conv}
\end{subfigure}
\caption{In \subref{fig:complex_shaped_sol}, the finite element solution of the problem with complex-shaped coarse mesh is shown, and in \subref{fig:complex_shaped_conv} is the convergence of said problem based on the number of layers the patches spread. Here, $A^0$ is a piecewise constant, random coefficient within $[0.01, 1]$, $A^1=1$, $B^1=1$, $f^0~=~1$ and $f^1=1$. We have $h=\frac{1}{32}$ and $H=\frac{1}{2}$.} \label{fig:complex_shaped_sol_conv}
\end{figure}
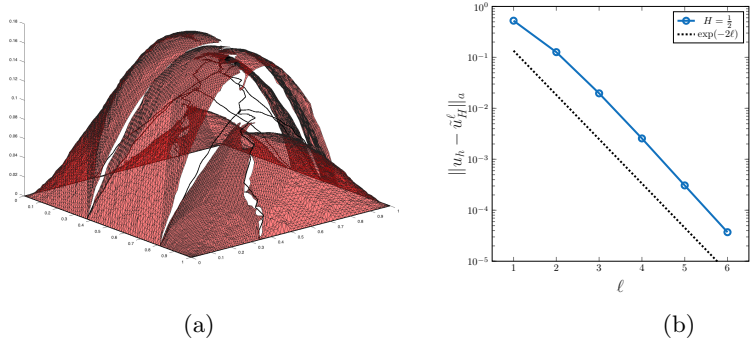

In the third numerical example, we consider a generalized coarse mesh consisting of complex-shaped elements; see \cref{sec:complexelements}. The domain and corresponding coarse mesh are shown in \cref{fig:complex_shaped_domain}\subref{fig:complex_shaped_domain_new}. In this experiment, $A^0$ is again chosen to be piecewise constant on the fine mesh $\mathcal{T}_h$, with random values in the range $[0.01, 1]$, while $A^1 = 1$, $B^1 = 1$, $f^0 = 1$, and $f^1 = 1$. The construction of the method for these complex-shaped coarse meshes follows the procedure described in \cref{sec:complexelements}.
For implementation purposes, we deviate from the stabilized localization strategy presented in \cref{sec:localization} and instead employ a simplified localization approach, where the global problems \cref{pbphiE} defining the basis functions are localized directly. Note that this method is affected by the previously described deterioration of errors if the oversampling parameter is not increased logarithmically as the mesh size decreases. The resulting solution and convergence behavior for $h = \tfrac{1}{32}$ and $H = \tfrac{1}{2}$, with respect to the oversampling parameter $\ell$, are shown in \cref{fig:complex_shaped_sol_conv}. We emphasize that for this choice of $f^0$ and $f^1$, the first term on the right-hand side of error estimate~\cref{eq:errestpracH1} vanishes and only the exponentially decaying localization error remains. The observed decay behavior is again consistent with \cref{thm:convergenceprot}.

\begin{figure}[h!] \centering
\includegraphics[scale=.25]{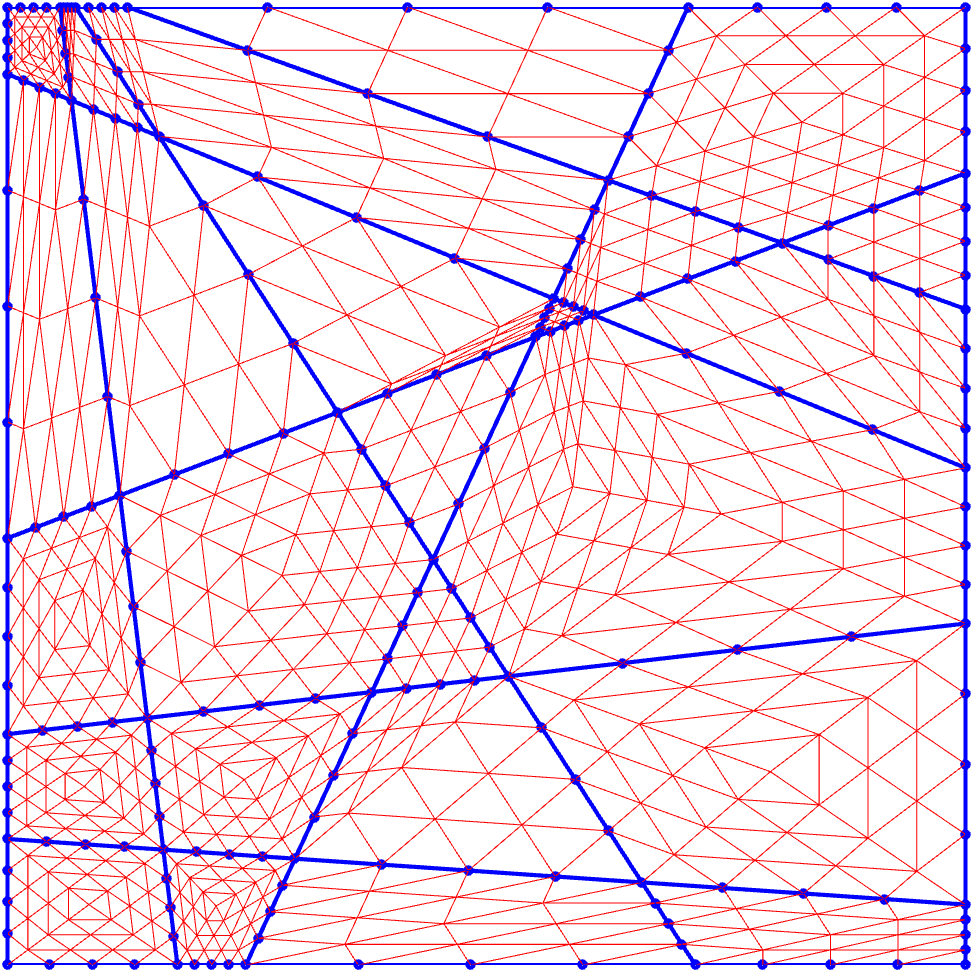}
\caption{The domain $\Omega$ used for the final numerical experiment. Here, the blue lines are the interfaces $\Omega^1$ and the coarse mesh corresponding to $H = \tfrac14$ is shown in red.}\label{fig:num_domain_8lines}
\end{figure}

\begin{figure}[h!]  \centering
\begin{subfigure}{.49 \textwidth} \centering
	\includegraphics[width=\textwidth]{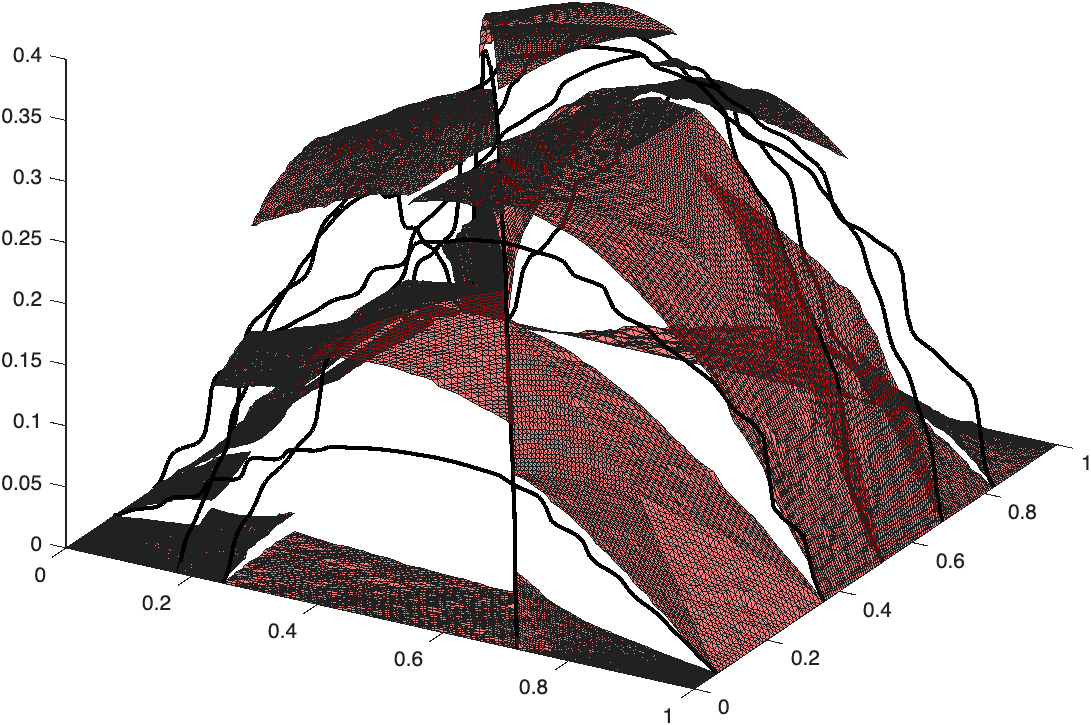}
	\caption{}\label{fig:sol_new_complicated1.1_h6}
\end{subfigure}
\begin{subfigure}{.49 \textwidth} \centering
        %
        %
        \definecolor{mycolor1}{rgb}{0.06600,0.44300,0.74500}%
        \definecolor{mycolor2}{rgb}{0.86600,0.32900,0.00000}%
        \definecolor{mycolor3}{rgb}{0.92900,0.69400,0.12500}%
        \definecolor{mycolor4}{rgb}{0.52100,0.08600,0.81900}%
        \definecolor{mycolor5}{rgb}{0.12941,0.12941,0.12941}%
        \begin{tikzpicture}[scale=.38]
        
        \begin{axis}[%
        width=3.5in,
        height=3.5in,
        at={(1.288in,0.657in)},
        scale only axis,
        xmode=log,
        xtick={0.125,0.25,0.5,1},
        xticklabels={{1/16},{1/8},{1/4},{1/2},{1}},
        xminorticks=true,
        xlabel style={font=\color{mycolor5}\huge, align=center},
        xlabel={$H$},
        ymode=log,
        yminorticks=true,
        ylabel style={font=\color{mycolor5}\huge},
        ylabel={$\|u_h - \tilde{u}^\ell_{H,h}\|_a$},
        axis background/.style={fill=white},
        legend style={at={(0.02,0.98)}, anchor=north west, legend cell align=left, align=left}
        ]
        \addplot [color=mycolor1, line width = 2, mark=o, mark size = 1, mark options={solid, mycolor1}]
          table[row sep=crcr]{%
        1	0.118827860813359\\
        0.5	0.048937681439517\\
        0.25	0.012740118625685\\
        0.125	0.003425259033119\\
        };
        \addlegendentry{$f^0 = \sin(\pi x) \sin(\pi y)$,$ f^1 = x + 2y$}
        
        \addplot [color=mycolor2, mark=o, line width = 2, mark size = 2, mark options={solid, mycolor2}]
          table[row sep=crcr]{%
        1	0.020761954519147\\
        0.5	0.012063336804951\\
        0.25	0.00725859769925\\
        0.125	0.002652453302586\\
        };
        \addlegendentry{$f^0 = f^1 = \sin(30\pi x)\sin(\pi y)$}
        
        \addplot [color=mycolor1, dotted, line width = 2]
          table[row sep=crcr]{%
        1	0.3\\
        0.5	0.075\\
        0.25	0.01875\\
        0.125	0.0046875\\
        };
        \addlegendentry{                $H^2$}
        
        \addplot [color=mycolor2, dotted, line width = 2]
          table[row sep=crcr]{%
        1	0.015\\
        0.5	0.0075\\
        0.25	0.00375\\
        0.125	0.001875\\
        };
        \addlegendentry{                $H$}
        
        \end{axis}
        \end{tikzpicture}%
	\caption{}\label{fig:convergence_new_complicated1.1_h6}
\end{subfigure}
\caption{In \subref{fig:sol_new_complicated1.1_h6}, the finite element solution is shown that corresponds to the problem with piecewise constant $A^0$ on the fine mesh, with values in the range $[0.01,~1]$. Here, $A^1 = \sin(30\pi x) \sin(30\pi y) + 1.1$ in order for the multiscale features to appear more clearly, whereas the other data is as in the previous examples. In \subref{fig:convergence_new_complicated1.1_h6}, we see the convergence with respect to $H$, following $\ell = \log_2(\frac{1}{H})$. The convergence for a smooth right hand side follows $H^2$ whereas using a more rapidly varying right hand side, we obtain a lower order of convergence.}\label{fig:new_complicated1.1_h6}
\end{figure}

In the final numerical example, we instead use the domain shown in Figure \ref{fig:num_domain_8lines}, utilising more interfaces for higher complexity. We let $A_0$ be piecewise constant on the fine mesh with random values in $[0.01,~1]$, $A^1=\sin(30 \pi x)\sin(30 \pi y)~+~1.1$, $B^1=1$, $f^0 = \sin(\pi x) \sin(\pi y)$ and $f^1= x + 2y$. Altering $A^1$ in this manner, we can more clearly see the multiscale effects of the solution on the interfaces, see \cref{fig:new_complicated1.1_h6} \subref{fig:sol_new_complicated1.1_h6}. \cref{fig:new_complicated1.1_h6} \subref{fig:convergence_new_complicated1.1_h6} shows the convergence with regards to $H$, letting $\ell = \log_2(1/H)$. As expected from \cref{thm:convergenceloc}, we obtain $H^2$-convergence using the smooth right hand side. For contrast, we also let $f^{0,1} = \sin(30 \pi x) \sin(\pi y)$, which is rapidly varying and not resolved by $H$ but is resolved by $h$. For this case, we get a lower order of convergence, which is expected according to \cref{thm:convergenceloc}.

\section*{Acknowledgments}

M.~Hauck acknowledges funding from the Deutsche Forschungsgemeinschaft\linebreak  (DFG, German Research Foundation) -- Project-ID 258734477 -- SFB 1173. Furthermore,  A.~M\aa lqvist acknowledges funding from the Swedish Research Council (VR) -- Project-ID 2023-03258\_VR.

\printbibliography

\end{document}